 \documentclass[preprint,3p,times,twocolumn]{elsarticle}
 
\usepackage{amssymb}
\usepackage[fleqn]{amsmath}
\DeclareMathOperator{\atantwo}{atan2}
\usepackage{booktabs}
\usepackage[skip=1pt]{subcaption}
\usepackage{siunitx}
\usepackage{bm}
\usepackage{lineno}

\usepackage{acro}

\NewAcroTemplate[list]{mythree}{% 
	\acronymsmapF{%
		\AcroAddRow{%
			\acrowrite{short}%
			& 
			&
			\acrowrite{long}%
			\tabularnewline
			\vspace{1mm}%
		}%
	}
	{\AcroRerun}%
	\acroheading
	\acropreamble
	\par\noindent
	\begin{tabular}{lll}
		\AcronymTable
	\end{tabular}
}
\DeclareAcronym{ABEP}{
	short = ABEP,
	long = air-breathing electric propulsion
}

\DeclareAcronym{ATLAS}{
	short = ATLAS,
	long = Advancing Technologies of Very Low-Altitude Satellites
}

\DeclareAcronym{CAM}{
	short = CAM,
	long = collision avoidance manoeuvre
}

\DeclareAcronym{CARA}{
	short = CARA,
	long = Conjunction Assessment Risk Analysis
}

\DeclareAcronym{CCSDS}{
	short = CCSDS,
	long = Consultative Committee for Space Data Systems
}

\DeclareAcronym{CDM}{
	short = CDM,
	long = conjunction data message
}

\DeclareAcronym{DIRECT}{
	short = DIRECT,
	long = Direct Shooting Method
}

\DeclareAcronym{COLA}{
	short = COLA,
	long = collision avoidance
}

\DeclareAcronym{DA}{
	short = DA,
	long = differential algebra
}

\DeclareAcronym{DM}{
	short = DM,
	long = descent metric
}

\DeclareAcronym{DSMC}{
	short = DSMC,
	long = Direct Simulation Monte Carlo
}

\DeclareAcronym{EIVE}{
	short = EIVE,
	long = Exploratory In-Orbit Verification of an E-band Link satellite
}

\DeclareAcronym{ESA}{
	short = ESA,
	long = European Space Agency
}

\DeclareAcronym{EUSST}{
	short = EUSST,
	long = European Union Space Surveillance and Tracking
}

\DeclareAcronym{FMF}{
	short = FMF,
	long = free-molecular flow
}

\DeclareAcronym{GEO}{
	short = GEO,
	long = Geostationary Orbit
}

\DeclareAcronym{GNSS}{
	short = GNSS,
	long = Global Navigation Satellite System
}

\DeclareAcronym{GOCE}{
	short = GOCE,
	long = Gravity Field and Steady-State Ocean Circulation Explorer
}

\DeclareAcronym{GSI}{
	short = GSI,
	long = gas-surface interaction
}

\DeclareAcronym{GSOC}{
	short = GSOC,
	long = German Space Operations Center
}

\DeclareAcronym{HBR}{
	short = HBR,
	long = hard-body radius
}

\DeclareAcronym{LEO}{
	short = LEO,
	long = low Earth orbit
}

\DeclareAcronym{MASTER}{
	short = MASTER,
	long = Meteoroid and Space Debris Terrestrial Environment Reference
}

\DeclareAcronym{MTCOLA}{
	short = MT-COLA,
	long = Minimum Lead Time Differential Algebra Method 
}

\DeclareAcronym{NASA}{
	short = NASA,
	long = National Aeronautics and Space Administration
}

\DeclareAcronym{NLP}{
	short = NLP,
	long = nonlinear program
}

\DeclareAcronym{OCP}{
	short = OCP,
	long = optimal control problem
}

\DeclareAcronym{ROAR}{
	short = ROAR,
	long = Rarefied Orbital Aerodynamics Research
}

\DeclareAcronym{SDMR}{
	short = SDMR,
	long = Space Debris Mitigation Requirements
}

\DeclareAcronym{SMA}{
	short = SMA,
	long = semi-major axis
}

\DeclareAcronym{SMD}{
	short = SMD,
	long = squared Mahalanobis distance
}

\DeclareAcronym{SOAR}{
	short = SOAR,
	long = Satellite for Orbital Aerodynamics Research
}

\DeclareAcronym{SRP}{
	short = SRP,
	long = solar radiation pressure
}

\DeclareAcronym{SSA}{
	short = SSA,
	long = space situational awareness
}

\DeclareAcronym{SSO}{
	short = SSO,
	long = sun-synchronous orbit
}

\DeclareAcronym{TCA}{
	short = TCA,
	long = time of closest approach
}

\DeclareAcronym{TPoC}{
	short = TPoC,
	long = total probability of collision
}

\DeclareAcronym{VLEO}{
	short = VLEO,
	long = very low Earth orbit
}

\acsetup{
	make-links = false,
	list/heading = section*,
	list/name={Abbreviations},
	list/template = mythree
}

\usepackage[unicode=true,
bookmarks=true,bookmarksnumbered=false,bookmarksopen=false,
breaklinks=false,pdfborder={0 0 1},backref=false,colorlinks=true]
{hyperref}

\usepackage[capitalise,nameinlink,noabbrev]{cleveref} %
\crefname{chapter}{Chapter}{Chapters}
\crefname{section}{Section}{Sections}
\crefname{subsection}{Section}{Sections}
\crefname{figure}{Fig.}{Figs.}
\crefname{table}{Table}{Tables}
\crefname{equation}{Eq.}{Eqs.}

\newlength{\figwidth}
\journal{Preprint submitted to Acta Astronautica}

\usepackage{xcolor}
\myfooter[C]{\normalsize\color{gray} Preprint submitted to Acta Astronautica - July 2026}

\begin{document}
	
\begin{frontmatter}
    \title{Aerodynamic Collision Avoidance Manoeuvres in Very Low Earth Orbit: Exploiting Lift for Short Warning Times and Reduced Orbital Decay}
	
	\author[aff1]{Paul J. Haufe}
	\author[aff1]{Fabrizio Turco}
	\author[aff2]{Friedrich Tuttas}
	\author[aff3]{Zeno Pavanello}
	\author[aff1]{Stefanos Fasoulas}
	\author[aff4]{Roberto Armellin} 
	
	\affiliation[aff1]{organization={Institute of Space Systems, University of Stuttgart},
		addressline={Pfaffenwaldring 29}, 
		city={Stuttgart},
		postcode={70569}, 
		country={Germany}}
	
	\affiliation[aff2]{organization={Institute of Flight Mechanics and Control, University of Stuttgart},
		addressline={Pfaffenwaldring 27}, 
		city={Stuttgart},
		postcode={70569}, 
		country={Germany}}

	\affiliation[aff3]{organization={Department of Aerospace Science and Technology (DAER), Politecnico di Milano},
		addressline={Via Giuseppe la Masa 34}, 
		city={Milan},
		postcode={20156}, 
		country={Italy}}
	
	\affiliation[aff4]{organization={Te Pūnaha Ātea - The Space Institute, The University of Auckland},
		addressline={20 Symonds Street}, 
		city={Auckland},
		postcode={1010}, 
		country={New Zealand}}
        
	%% Abstract
	\begin{abstract}
		Very Low Earth Orbits (\acs{VLEO}, below \qty{450}{\kilo\metre} altitude) offer promising mission capabilities for future spacecraft, but increasing traffic and the presence of orbital debris create collision risks that pose a particular challenge for satellites without conventional propulsion. Aerodynamic drag, resulting from interaction with the residual atmosphere, has been shown effective for performing \acp{CAM}. However, short warning times remain challenging, and current methods do not exploit aerodynamic lift forces that enable out-of-track manoeuvres, which could help enhance control authority and mitigate drag-induced orbital decay. This work addresses both gaps by developing and evaluating two complementary optimisation methods: a heuristic approach using reachable set analysis to determine minimum lead times, and a direct shooting method to limit the increase in orbital decay due to the manoeuvre. Analyses using data derived from real conjunctions show that achieving a sufficient reduction in collision risk in the considered scenarios is feasible with lead times lower than one orbital period. Results also show that lift-based manoeuvres outperform drag-only approaches, with modest lead time reductions and more pronounced benefits in limiting orbital decay. Performance sensitivity to variations in control authority due to the type of gas-surface-interaction, altitude, and the encounter geometry is characterised. The methods and results presented highlight the potential of aerodynamic \acp{CAM} for short warning time scenarios, contributing to sustainable operations in the increasingly congested \ac{VLEO} environment.
	\end{abstract}
	
	
	%% Keywords
	
\end{frontmatter}
% -------------------------------------------------------------------------------------%
% -------------------------------------------------------------------------------------%
% ----------------------------------- Introduction ------------------------------------ %
% -------------------------------------------------------------------------------------%
% -------------------------------------------------------------------------------------%
\acresetall
\section{Introduction}
%\linenumbers
The near‑Earth environment has become progressively congested by space traffic and debris. Collisions pose a lethal risk to operational spacecraft and require actions in form of \acp{CAM}. Increasing launch traffic and large constellations accelerate the trend despite upcoming mitigation efforts and guidelines \cite{esa2023}. Very Low Earth Orbits (\acs{VLEO}, below \qty{450}{\kilo\metre}) offer operational benefits with higher resolution in Earth observation, lower communication latency and improved link budgets. Although the faster natural decay aids the removal of debris, the rising spatial density due to constellations also affects the upper end of \ac{VLEO} as captured by the \ac{ESA}'s \textit{MASTER} population \cite{spacedebrisoffice2025}. Consequently, collision risks and frequent conjunctions remain operational challenges even at very low altitudes.

Aerodynamic manoeuvres exploit forces from interaction with the residual atmosphere as a propellant‑free means to change a satellite’s trajectory. Drag‑based \acp{CAM} act primarily in-track and have been shown effective when planned several orbits ahead \cite{bella2025,turco2023}. However, short warning times and increased orbital decay due to the manoeuvre present open challenges that could be addressed by using lift forces acting perpendicular to drag. Advances in shape optimisation and surface materials promise an increase in achievable aerodynamic lift \cite{hild2022}, enabling out‑of‑track control that is not possible with drag alone. This is especially valuable for late, just‑in‑time \acp{CAM} performed only a few orbits before the projected conjunction.  Fresher tracking data reduces propagated uncertainty and helps avoid unnecessary manoeuvres \cite{nasa2023}.

This work develops methods to design aerodynamic \acp{CAM} that include lift and quantifies their benefits over drag‑only manoeuvres. Two complementary solution methods are presented: a reachable set heuristic that is  based on a backward sweep approach and identifies minimum feasible lead times, and a direct shooting approach that minimises the inevitable increase in orbital decay due to the manoeuvre. Problems of \ac{CAM} design are often solved using direct and indirect methods \cite{pavanello2025a}. However, since the dynamics as well as the constraint on collision risk are nonlinear and, in general, non-convex, convergence to the global optimum is not guaranteed. The developed heuristic method builds upon recent work on time-optimal \acp{CAM} \cite{dellelce2024, pavanello2026a} and extends them to manoeuvers with constrained control direction. It serves as an initial guess for further refinement using direct methods for minimising decay. The proposed tools operate directly on collision risk and track changes in \ac{TCA} due to the manoeuvre. Performance is evaluated across representative conjunction geometries, altitudes and surface properties.

The paper is organised as follows: \cref{sec:problemDefinition} introduces relevant concepts for the \ac{CAM} design problem and aerodynamic manoeuvres. \cref{sec:methods} presents the developed approach and implementation, highlighting both the chosen solution methods to the minimum lead time and the minimum decay problems. The results and a discussion of lift-based \acp{CAM} is found in \cref{sec:results} and \cref{sec:conclusions} is a short summary and conclusion.

% -------------------------------------------------------------------------------------%
% -------------------------------------------------------------------------------------%
% ---------------------------------- Fundamentals ----------------------------------- %
% -------------------------------------------------------------------------------------%
% -------------------------------------------------------------------------------------%

\section{Aerodynamic Collision Avoidance}\label{sec:problemDefinition}
\begin{figure}[tb]
	\centering
	\includegraphics[width=\figwidth]{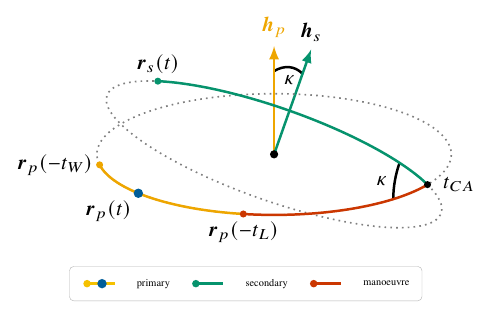}
	\caption[Illustration of the relevant parameters in a conjunction]{Conjunction setup with primary ($p$) and secondary spacecraft ($s$), warning time $t_W$, manoeuvre lead time $t_L$, time of closest approach $t_{CA}$ and mutual inclination $\kappa$.}
	\label{fig:ConjScenario}
\end{figure}

In \ac{CAM}-design, it is common to make the distinction between the primary maneuverable spacecraft ($p$) and one or multiple secondary objects ($s$), which are considered uncontrolled debris. For this work, we focus on a single conjunction with a single secondary object. \Cref{fig:ConjScenario} displays the conjunction setup and terminology used. The orbits of both objects are shown with the mutual inclination $\kappa$ as the angle between the angular momentum vectors $\bm h_p$ and $\bm h_s$. The \ac{TCA}, or $t_{CA}$, is the point in time at which the distance between the two objects is minimum. We consider short warning time scenarios, where the earliest time that a \ac{CAM} may begin ($t_W$) is on the order of a few orbits. The actual start time of the \ac{CAM} is the manoeuvre lead time $t_L$, and in the case of minimum-time collision avoidance, the manoeuvre lasts continuously until $t_{CA}$. Note that $t_W$ and $t_L$ are time spans, so their value is positive, while the time variable $t$ is set to zero at $t_{CA}$ and is negative prior, so 
\begin{equation}
	t\in\left[-t_W,0\right]\subset\mathbb{R}_{\le0}~.
\end{equation}

\subsection{Dynamics and B-plane}
We consider satellites in \ac{VLEO} and describe motion using the restricted two-body problem. The state vector of the primary object $\bm x_{p}$ is
\begin{equation}
	\begin{aligned}
		\bm x_{p}&= \begin{pmatrix} \bm r_{p} \\ \bm v_{p} \end{pmatrix}\in\mathbb{R}^{6}, 
	\end{aligned}
\end{equation}
with cartesian position $\bm r_{p}$ and velocity $\bm v_p$. The dynamics of the restricted two-body problem in the Earth-Centered Inertial frame~(I) and under aerodynamic perturbation is the following system of ordinary differential equations:
\begin{equation}\label{eq:sysOde}
	\dot{\bm{x}}_{p}(t) =
	\begin{pmatrix}
		\dot{\bm{r}}_{p}(t) \\
		\dot{\bm{v}}_{p}(t)
	\end{pmatrix}
	=
	\begin{pmatrix}
		\bm{v}_{p}(t) \\[2mm]
		-\mu_E \dfrac{\bm{r}_{p}(t)}{\|\bm{r}_{p}(t)\|^3} + \bm{f}_\textrm{aero}^{I}(t)
	\end{pmatrix},
\end{equation}
where $\mu_E$ is the gravitational parameter of Earth. An additional acceleration, known as the specific aerodynamic force $\bm{f}_\textrm{aero}^I$, results from the interaction between the residual atmosphere and the spacecraft's surface. This force is used to control the trajectory in aerodynamic manoeuvres. The secondary object's state $\bm x_{s}$ is defined analogously, but remains uncontrolled. 

Secondary perturbations due to contributions from higher-order harmonics of Earth’s gravitational potential and solar radiation pressure are neglected, since the goal of \ac{CAM} design is to compute necessary deviations from the nominal trajectory, not to predict absolute orbital motion accurately. Therefore, the influence of other perturbations only needs to be considered if the change in trajectory and the conjunction event achieved by the \ac{CAM} differs meaningfully in the presence of these perturbations. An analysis of Hernando-Ayuso and Bombardelli has shown that, in comparison with a high-fidelity propagator, contributions from third-body gravity and higher-order Earth gravitational harmonics are negligible for \ac{CAM} design in low Earth orbits and with a time horizon of only few orbits \cite{hernando-ayuso2021}. 

\begin{figure}[tb]
	\centering
	\includegraphics[width=\figwidth]{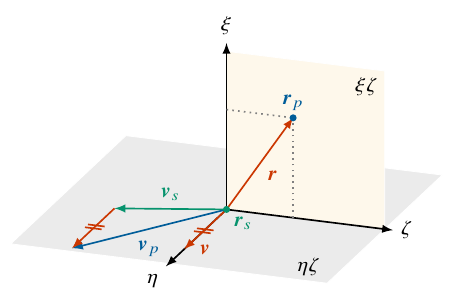}
	\caption[Conjunction geometry and B-plane coordinate system definition.]{Conjunction geometry and B-plane coordinate system definition with relative velocity $\bm{v}_{CA}$ and relative position $\bm{r}_{CA}$ at \ac{TCA} (adapted from \cite{armellin2021}). }
	\label{fig:CosBplane}
\end{figure}

The relative position and relative velocity of the two objects at \ac{TCA} are
\begin{equation}
	\begin{aligned}
		\bm r\,(t_{CA})&=\bm r_{p}(t_{CA})-\bm r_{s}(t_{CA}) \\
		\bm v\,(t_{CA}) &=\bm v_{p}(t_{CA})-\bm v_{s}(t_{CA})~.
	\end{aligned}
	\label{eq:relStates}
\end{equation}
At TCA, the distance between both objects is minimal and therefore 
\begin{equation}
	\bm r\,(t_{CA}) \cdot \bm v\,(t_{CA})=0~.
	\label{eq:orthoCon}
\end{equation} 
A useful coordinate system is constructed based on the relative velocity components at \ac{TCA}, with the origin at the centre of the secondary object. The axes of this B-plane system (see \cref{fig:CosBplane}) are defined in accordance with the literature \cite{hernando-ayuso2021,armellin2021} as
\begin{equation}
	\begin{aligned}	 
		\bm{\hat u}_{\xi}&= \frac{\bm v_{s}\times \bm v_{p}}{||\bm v_{s}\times \bm v_{p}||} \\
		\bm{\hat u}_{\eta}&= \frac{\bm v_{p}-\bm v_{s}}{||\bm v_{p}-\bm v_{s}||} \\
		\bm{\hat u}_{\zeta}&= \bm{\hat u}_{\xi}\times \bm{\hat u}_{\eta} ~,
	\end{aligned}
\end{equation}
so $\bm{\hat u}_{\eta}$ points along the relative velocity direction and $\bm{\hat u}_{\xi},~\bm{\hat u}_{\zeta}$ are perpendicular to it. Since the relative position in the \mbox{$\bm{\hat u}_{\eta}$-direction} at TCA is zero by definition, the conjunction can be expressed using the reduced B-plane coordinates 
\begin{equation}\label{eq:redBplaneCoord}
\bm r^{\star}(t_{CA})=(\xi,\zeta)~. 
\end{equation}
The relative position at \ac{TCA} is the basis for computing the risk of both objects colliding.

\subsection{Probability of Collision}
During conjunction assessment and \ac{CAM} design, spacecraft are often conservatively modelled as spheres identified by the hard-body radius $R_{HB}$. If the relative position at \ac{TCA} $\bm{r}_{CA}$ is smaller than the combined $R_{HB}$, a direct hit occurs. Since the state uncertainty of orbiting spacecraft is high, it is common to assess collision risk statistically rather than deterministically. The probability of collision, or $P_C$, serves as a metric for spacecraft operators to form decisions about whether a \ac{CAM} shall be performed and how effective it is. 

We consider a short-term encounter, characterised by a high relative velocity. The encounter has a short duration, therefore it is common to approximate the relative motion of the objects as rectilinear and disregard velocity uncertainties \cite{chan2004}. Moreover, the assumption is that the position uncertainty follows a Gaussian distribution and remains constant for the duration of the encounter. Therefore, the relative position at \ac{TCA} follows a Gaussian distribution with mean $\bm{\mu}_{CA}$ and combined covariance of both objects $\mathbf{\Sigma}_{CA}$ as the sum of the individual covariances in a common frame. Within this work, both are given in reduced B-plane coordinates since we only consider the instantaneous probability of collision at \ac{TCA}: the $\eta$ component of $\bm{\mu}_{CA}$ is zero by definition.

Within this work, a \ac{CAM} is deemed successful if it reduces $P_C$ below a safe threshold. For satellites in low Earth orbit that comply with the \ac{ESA} requirements, a manoeuvre is required above $10^{-4}$ and must reduce the immediate risk to below $10^{-6}$ \cite{esa2023}. Several models to compute $P_C$ exist in the literature. For \ac{CAM} optimisation, a useful analytical approximation by Alfriend and Akella \cite{alfriend1999} assumes a constant probability density over $R_{HB}$ 
\begin{equation}
	P_C = \frac{R_{HB}^2}{2\sqrt{\det\,\mathbf{\Sigma}_{CA}}} \exp\left(-\frac{d^2}{2}\right)~,
\end{equation}
where the \ac{SMD} is defined as
\begin{equation}\label{eq:SMD}
	d^2 = \bm{\mu}_{CA}^{T}\,(\mathbf{\Sigma}_{CA})^{-1}\,\bm{\mu}_{CA}\quad\in\mathbb{R}_{+}~.
\end{equation}
This approximation deviates from higher-fidelity methods for large risk values. Therefore, to validate the success of a CAM, we use Chan's method of equivalent cross-sections to compute a higher-fidelity representation of $P_C$ (see \cite{chan2018}). The metric of \ac{SMD} is still useful since it can serve as a proxy for the collision risk gradient during optimisation through the closed-form solution:
\begin{equation}
	\nabla_{\bm{\mu}}\, d^2 = 2\,(\mathbf{\Sigma}_{CA})^{-1}\,\bm{\mu}_{CA}~.
\end{equation}
As Chan's method enables the computation of iso-probability levels on the B-plane via inversion, \ac{SMD} and $P_C$ can be related \cite{devittori2026}. The initial conditions for \ac{CAM} design, namely the state vectors, covariance matrices and hard body radii of both objects at \ac{TCA}, are taken from \acp{CDM} issued by monitoring entities. For this work, we base our analyses on data of selected representative conjunctions from the \textit{ESA Collision Avoidance Challenge}, listed in \cref{tab:mtColaScenarios} \cite{esa2019}.

\begin{figure}[tb]
	\centering
	\includegraphics[width=0.95\figwidth]{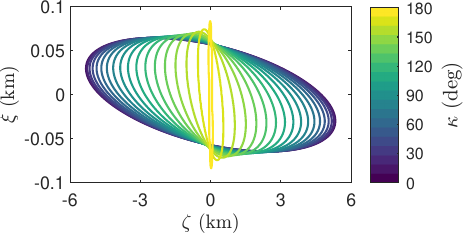}
	\label{fig:isoProb77}	
	\caption[Iso-probability ellipses for chosen conjunction events and a variation in mutual inclination]{Iso-probability ellipses for $P_C=10^{-6}$ and variation in mutual inclination $\kappa$ for encounter $I_\textrm{CDM}=77$ in \cref{tab:mtColaScenarios}.}
	\label{fig:isoProb}
	\vspace{-2mm}
\end{figure}

Given the covariances of the two objects at \ac{TCA}, we can visualise iso-probability levels that form ellipses on the B-plane for the chosen $P_C$ approximation. \cref{fig:isoProb} shows the $P_C=10^{-6}$ level for variation in mutual inclination $\kappa$ for encounter $I_\textrm{CDM}=77$ in \cref{tab:mtColaScenarios}. For moderate mutual inclinations, iso-probability ellipses are often elongated along $\zeta$. This axis corresponds to the along-track component of the positional covariance, which is typically largest, especially in low orbits under the influence of atmospheric drag \cite{reihs2017}. They may additionally be slanted depending on the individual covariance matrices of both objects and $\kappa$. The tilt of the combined covariance is relevant since it influences the B-plane position change needed to reach safe thresholds. As described by Hernando-Ayuso and Bombardelli in the context of tangential low-thrust encounters, some \ac{CAM} scenarios change drastically due to this tilt \cite{hernando-ayuso2021}. As can be seen in \cref{fig:isoProb}, the ellipses widen for growing $\kappa$ along $\xi$ and contract along $\zeta$. Head-on collisions ($\kappa=\ang{180}$) are special, since the velocities are exactly parallel and the B-plane coordinate system is undefined about the $\eta$ axis. Then the choice of $\xi$ and $\zeta$ orientation is arbitrary and the along-track components of the individual covariances are oriented in the direction of $\eta$, thus have no effect on the iso-probability ellipses. 

\subsection{Satellite Aerodynamics}
\begin{figure}[tb]
    \centering
    \includegraphics[clip, trim={2mm 3mm 2mm 2mm},width=\figwidth]{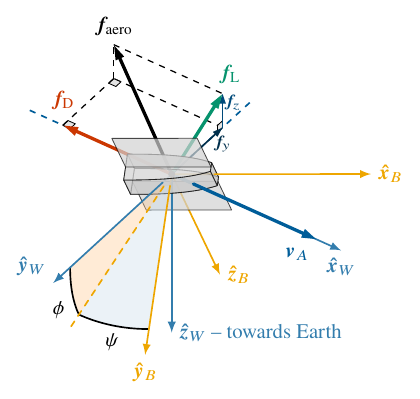}
	\caption{Wind (W) and body (B) coordinate systems with the rotation sequence ($\phi\rightarrow\psi$) that relates them, and illustration of lift $\bm{f}_\mathrm{L}$, drag $\bm{f}_\mathrm{D}$, and the resulting aerodynamic force $\bm{f}_\mathrm{aero}$.}
    \label{fig:CosAeroBody}
\end{figure}

Forces resulting from aerodynamic interaction with the residual atmosphere are the dominant perturbation for objects in \ac{VLEO}. Satellites can exploit this for aerodynamic manoeuvres by varying the incident angle and area of exposed surfaces through attitude changes or movable panels. Aerodynamic coefficients are a common approach for relating the satellite's attitude to aerodynamic forces while keeping the reference area fixed. For this, we define a wind (W) reference frame, as shown in \cref{fig:CosAeroBody}, that has the $\bm{\hat x}_{W}$ axis coincide with the velocity of the satellite relative to the atmosphere $\bm v_A$. The axis $\bm{\hat z}_{W}$ is defined via the normalised cross product $ \bm h_p \times \bm{\hat x}_{W}$ with the orbit's angular momentum vector $\bm h_p$, and $\bm{\hat y}_{W}$ completes the right-handed coordinate system. In circular orbits with no wind, $\bm{\hat z}_{W}$ lies in the orbital plane pointing towards Earth, being anti-parallel to the position vector. The body frame is obtained from the wind frame by rotating about $\bm{\hat x}_W$ by $\phi$, followed by rotating about the new $\bm{\hat z}_W$ by $\psi$. The specific aerodynamic force $\bm f_\textrm{aero}^{W}$, is calculated using the force coefficients $\bm{C}_{A}^W=[C_x, C_y, C_z]^T$ in the wind frame with
\begin{equation}\label{eq:aeroForce}
	\bm{f}_\textrm{aero}^{W}= \frac{\rho ~A_{ref}}{2~m}~||\bm{v}_{A}||^{2}~ \bm{C}_{A}^W~,
\end{equation}
where $\rho$ is the mass density of the atmosphere, $A_{ref}$ is a reference area matched to the coefficients and $m$ the satellite mass.

High-fidelity density models account for the influence of space weather, current latitude, and seasonal variations. However, this study aims to assess the qualitative use of lift and drag in aerodynamic \acp{CAM}. Therefore, we use the exponential atmosphere model, which captures the influence of variations in density due to the aforementioned factors by varying the altitude only. Effects of the co-rotation of Earth's atmosphere and winds are not considered in order to limit the number of parameters affecting the subsequent analyses. Since wind effects are not considered in this work, the wind reference frame is fixed with respect to the in-track satellite-centered coordinate system (NTW), with velocity direction \mbox{$\bm{\hat x}_{W}=\bm{\hat T}$}, negative angular momentum vector  \mbox{$\bm{\hat y}_{W}=-\bm{\hat W}$} and \mbox{$\bm{\hat z}_{W}=-\bm{\hat N}$}. 

It is useful to distinguish between drag and lift forces. Drag $\bm f_D$ always acts in negative $\bm{\hat x}_W$ since $C_x^W\in\mathbb{R}_{-}$, while lift $\bm f_L$ remains in the $\bm{\hat y}_{W}\bm{\hat z}_{W}$-plane. The lift coefficient $C_L\in\mathbb{R}_{+}$ is the Euclidean norm of the sum of $C_y$ and $C_z$, which act in $\bm{\hat y}_W$ and $\bm{\hat z}_W$ directions, respectively. The lift and drag magnitudes are therefore determined by the value of $\psi$, while the produced lift vector rotates around $\bm{\hat x}_W$ via $\phi$. An important metric for aerodynamic manoeuvres is the lift-to-drag ratio with 
\begin{equation}
	L/D= \frac{f_L}{f_D} = \frac{C_L}{C_D} = \frac{\sqrt{C_y^2+C_z^2}}{-C_x}~.
\end{equation}
New satellites designed for \ac{VLEO} promise increased control authority through improved $L/D$ ratios and offer the ability to perform prolonged manoeuvres by staying in off-nominal attitudes via active aerodynamic attitude control and trimming. The concept remains relatively new, with only a few missions providing limited in-orbit data that could be used to draw conclusions on realistic lift magnitudes (e.g. \cite{crisp2022a}).
Aerodynamic coefficients for realistic satellite geometries in \ac{VLEO} are commonly estimated using numerical methods. High-fidelity approaches, such as \ac{DSMC}, aim to capture particle interactions accurately. While highly detailed \ac{DSMC} analysis is often too computationally expensive for design-space exploration, an alternative class of methods is panel-based. The spacecraft geometry is discretised into flat panels to compute particle incident angles as well as shadowing. For tools like \textit{ADBSat}, this enables rapid computation of aerodynamic coefficients for quasi-realistic geometries and a variety of flow directions \cite{sinpetru2022a}. 

Regardless of their computational approach, all methods rely on \ac{GSI} models that capture the extent to which the momentum from an incident particle is transferred to the spacecraft surface. Parameters such as surface properties, temperature, and incident angle influence the direction and velocity of the reflected or re-emitted particle, and the resulting reaction force on the vehicle \cite{livadiotti2020}. A simple model that is available both in \textit{ADBSat} and \textit{PICLas} is that of \textit{Maxwell} \cite{maxwell1878}. It models the particle interaction using a single coefficient $\sigma$ that describes the fraction of diffuse particle re-emission. In short, $\sigma=1$ corresponds to fully diffuse re-emission and $\sigma=0$ is the case for specular reflection \cite{livadiotti2020}. The specular property offers the highest theoretical $L/D$ ratios, but no durable materials with sustained specular behaviour have been demonstrated to date. Atomic oxygen in the atmosphere reacts with materials, thereby degrading the reflective property of all exposed surfaces \cite{moe2005}.

The spacecraft geometry plays a vital role in increasing aerodynamic control authority. Various studies with the aim to reduce drag in the nominal flight configuration or to increase the L/D ratio have been conducted \cite{hild2022,sinpetru2022,agez2025}. Marianowski et al. optimised the 3-unit CubeSat shown in \cref{fig:CosAeroBody} for differential lift and drag control \cite{marianowski2024}. Turco et al. performed additional analyses for this geometry and varying $\sigma$ to study aerodynamic manoeuvres in inclination changes with minimal decay in the \ac{SMA} \cite{turco2025}. Next, these coefficients are used to construct a tangible representation of the aerodynamic control authority.

\subsection{Aerodynamic Control Authority}
\begin{figure}[tb]
	\centering	
	\includegraphics[width=0.85\figwidth]{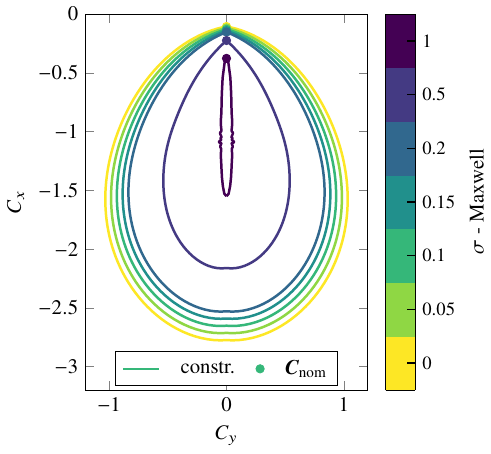}
	\caption[Section of the fitted polynomial control constraints]{Section of the fitted polynomial control constraints.}
	\label{fig:aeroConstr}		
\end{figure}

\begin{figure}[tb]
	\centering
	\hspace{0cm}	
	\includegraphics[clip,trim={3mm 3mm 3mm 3mm},width=0.95\figwidth]{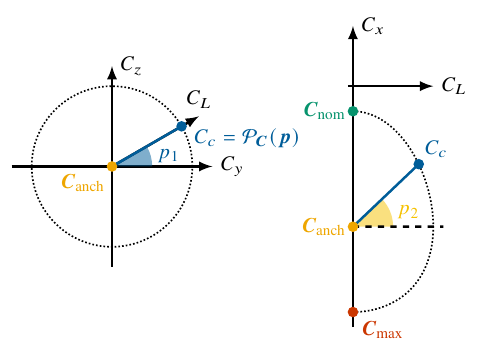}
	\caption[Illustration of the pseudo-controls and their use in the aerodynamic control constraint]{Illustration of the pseudo-controls $\left(p_1,p_2\right)$ and their use in the aerodynamic control constraint.}
	\label{fig:aeroConstrPseudoCtr}		
\end{figure}

For use in the developed optimisation methods, functions that capture the complex constraints on aerodynamic control acceleration are constructed. In this study, we modulate $\bm C_A$  via changes in attitude by $\phi$ and $\psi$-rotations. As the relative velocity always remains in the body's symmetry plane, sideslip is zero, and only the coefficients for variation in $\psi$ are required. Aerodynamic coefficients for a realistic geometry were derived by using the Maxwell \ac{GSI} model and \textit{ADBSat}. An extended coefficient grid to that in \cite{turco2025} was generated for 
\begin{itemize}
\item $\sigma\in\{0, 0.05, 0.10, 0.15, 0.2, 0.5, 1\}$, and
\item $\psi\in[\ang{-90},\ang{90}]$ in $\ang{1}$ increments.
\end{itemize}
The reference area used is ${A_{ref}=\qty{0.13}{\meter\squared}}$ and the satellite's mass is $\qty{5}{\kilogram}$. The resulting coefficient curves are shown in  \Cref{fig:aeroConstr}, with nominal control $\bm C_\textrm{nom}$ in the minimum drag configuration and the maximum drag is reached for $\psi=\pm\ang{90}$. An anchor point $\bm C_\textrm{anch}$ is chosen as the mean of the coefficient data, so the shifted point is
\begin{equation}
    \bm C_A' = \bm C_A - \bm C_\textrm{anch}~.
\end{equation}
The control constraint is based on a parametrisation in spherical coordinates with pseudo-controls $\bm p=\left(p_1,p_2\right)^T$ as shown in \cref{fig:aeroConstrPseudoCtr}. Given the aerodynamic coefficient data $\bm C'_A$, spherical coordinates are constructed as
\begin{equation}
	\begin{array}{c c r l l }
		&p_1 &=& \atantwo(C'_z, C'_y)~,		\\[1mm]
		&p_2 &=& \atantwo(C'_x, C'_{L})~,\\[2mm]
		\text{with}&	C'_{L} &=& \sqrt{(C'_{y})^{2} + (C'_{z})^{2}}& \text{, the lift magnitude.}\\		
	\end{array}
\end{equation}
A bivariate polynomial $C_{c}$ that maps pseudo-controls to aerodynamic coefficient magnitude is obtained from the point cloud via least-squares fitting. The resulting map 
\begin{equation}
	C_{c}=\mathcal{P}_{\bm C}(\bm p)~\quad\textrm{with} \quad p_1\in[-\pi,\pi],\quad p_2\in\left[-\frac{\pi}{2},\frac{\pi}{2}\right].
	\label{eq:constraint}
\end{equation}
yields the aerodynamic coefficient vector for given pseudo-controls as 
\begin{equation}
	\bm{C}_A = 
	\left(\begin{array}{@{}l@{}}
		C_{\textrm{anch},x} + C_\textrm{c} \sin(p_2) \\
		C_{\textrm{anch},y} + C_\textrm{c} \cos(p_1) \cos(p_2) \\
		C_{\textrm{anch},z} + C_\textrm{c} \sin(p_1) \cos(p_2)
	\end{array}\right)~.
\end{equation}
Together with \cref{eq:aeroForce}, this completes the mapping from pseudo-control vector $\bm p$ to the aerodynamic force vector $\bm{f}_\textrm{aero}$.  The full chain of operations is captured as the composed control mapping function $\mathcal{C}_{f,\sigma}$ for each $\sigma$ as
\begin{equation}\label{eq:pseudoToAccel}
	\bm f_{\mathrm{aero},\sigma}^W= \mathcal{C}_ {f,\sigma}(\bm p).
\end{equation}
Polynomials as constraints offer benefits in continuity, differentiability and computational efficiency compared to querying the raw point cloud data directly during optimisation. The fully diffuse case is excluded from numerical analysis because the corresponding aerodynamic coefficient surface becomes highly elongated along the $C_x$-axis, resulting in a near-linear relation with very small lift components. This shape is inadequately represented by the angular pseudo-control parametrisation and leads to numerical ill-conditioning during the polynomial fitting, as seen in the oscillations for $\sigma = 1$ in \cref{fig:aeroConstr}. Alternative parametrisations, such as ellipsoidal or Fourier-based constraints, were tested but either overestimated $L/D$ in the low drag region or were incompatible with the developed routines and are therefore not implemented at this point. Additionally, the fully diffuse case is predicted to inhibit the lowest performance gains of lift due to the low $L/D$. Consequently, the analysis focuses on the remaining values of $\sigma$ that yield stable and well-behaved polynomial mappings. Excluding the fully diffuse case, the fits achieve a maximum relative magnitude error of $\le\qty{2.8}{\percent}$, with the largest deviations at high $\sigma$. 

In practical mission scenarios, a minimum nominal drag configuration, as shown in \cref{fig:aeroConstrPseudoCtr}, is most favourable to prolong the mission and is therefore taken as the basis for this work. For this case, the produced control accelerations are always in the negative flight direction with some cross-track component due to lift. The approach can be adapted to other nominal drag configurations, spacecraft geometries, and non-axisymmetric coefficient point clouds.
In the analysis of aerodynamic \ac{CAM}s, three cases are discussed and compared:
\begin{itemize}
	\item \textbf{Drag-only}: control is restricted to lie on the $C_x$-axis. 
	\item \textbf{Lift-based}: in addition to the induced drag, lift forces are generated, and their magnitude depends on the \ac{GSI} property.
	\item \textbf{Propulsive}: accelerations corresponding to the maximum drag coefficients $\bm C_\mathrm{max}$ can be generated that are not constrained in direction.
\end{itemize}
For the lift-based case, the constraint is egg-shaped and axisymmetric around $\bm{\hat x}_b$, while in the drag-only case, it becomes a line, and for propulsive manoeuvres a sphere centred in the nominal drag configuration. The maximum drag configuration $\bm C_\textrm{max}$ is reached for $p_2 = -\frac{\pi}{2}$. Using this parametrisation, $p_1$ is interchangeable with $\phi$, while the relation of $p_2$ to $\psi$ is captured through a polynomial map. 

% -------------------------------------------------------------------------------------%
% -------------------------------------------------------------------------------------%
% -------------------------------------- Method -------------------------------------- %
% -------------------------------------------------------------------------------------%
% -------------------------------------------------------------------------------------%

\section{Manoeuvre Optimisation Approach}\label{sec:methods}
Designing \acp{CAM} with minimum lead time and decay constitutes an optimal control problem, which we solve numerically by formulating it as a \ac{NLP}. While a multitude of methods exist for impulsive, low-thrust and even drag-based \acp{CAM} (see \cite{pavanello2025a}), solutions for the minimum lead time problem and the incorporation of aerodynamic lift remain unexplored. This section presents our novel methods and algorithms for these problems under aerodynamic constraints. We first construct polynomial maps of the conjunction dynamics to simplify problem setup. Next, we define the minimum lead time \ac{NLP} and describe the developed heuristic method for solving it, which includes a reachable set approximation that builds intuition and serves as a basis for the solution. Finally, we address the minimum decay problem using a direct shooting approach.

\subsection{Construction of Polynomial Maps}\label{sec:polynomialMaps}
First, we construct polynomial maps that capture the dynamics of the conjunction and account for the critical effect that control actions have on \ac{TCA}. This approach leverages \ac{DA} to efficiently compute high-order sensitivities without repeated numerical integration.

In astrodynamics, \ac{DA} is a computational framework that enables the automatic computation of high-order Taylor series expansions of the flow of dynamical systems up to a truncation order $k$. Instead of repeated numerical integrations, \ac{DA} can be used to construct polynomial maps that approximate how system states evolve with respect to initial conditions, parameters, and control inputs. The main advantage is efficiency: once the expansion is computed, system behaviour under perturbations can be analysed through simple function evaluations rather than expensive numerical propagation \cite{armellin2010}. \ac{DA} also enables direct computation of sensitivities, reuse of expansions across different scenarios, and effective treatment of problems where small forces accumulate over long time spans, such as low-thrust trajectory design \cite{pavanello2025b}. When solving \acp{NLP} for \acp{CAM}, the polynomial representation helps to efficiently capture how control variations propagate through the nonlinear dynamics for each discretised arc and affect the final collision risk. 

The interval between the warning time and \ac{TCA}, $[t_0,t_{CA}]$, is discretised into $N$ control arcs with discrete time nodes $\{t_n\}_{n=0}^{N}$. For notation purposes and to stay compatible with literature on propulsive \acp{CAM}, for which the developed tools can be used too, the control acceleration $\bm u$ is used interchangeably with the aerodynamic specific control force $\bm f_{\mathrm{aero}}$. Variations about the nominal trajectory of the primary $\bm{\bar x}_{p}$ and nominal control $\bm{\bar u}$ (i.e., minimum drag configuration) are introduced at each time step $t_n$ as
\begin{equation}
	\begin{aligned}
		\bm x_{p}(t_{n})&=\bm{\bar x}_{p}(t_{n})+\delta \bm x_{p}(t_{n}) \\ 
		\bm u(t_{n})&=\bm{\bar u}(t_{n}) + \delta \bm u(t_{n}),
	\end{aligned}
\end{equation}
and the system dynamics are expanded using \ac{DA}-based numerical integration up to truncation order $k$ (see \cite{armellin2010}). This yields Taylor polynomial maps of the flow at each step $t_{n}$ about the expansion point of the nominal trajectory $\bm{\bar x}_{p}$ and control $\bm{\bar u}$. The resulting polynomial map of the primary's state $\mathcal{T}_{\bm x_{p}}$ solely depends on the variations around the expansion point as  
\begin{equation}\label{eq:xpMap}
	\bm x_{p}(t_{n})=\mathcal{T}_{\bm x_{p}}\big(\delta \bm x_{p}(t_{n-1}),~\delta \bm u(t_{n-1})\big)~.
\end{equation}
If the constant part of the polynomial is omitted, the state and control deviation in a previous step relates to the deviation of the next state as
\begin{equation}\label{eq:dxpMap}
	\delta \bm x_{p}(t_{n})=\mathcal{T}_{\delta \bm x_{p}}\big(\delta \bm x_{p}(t_{n-1}),~\delta \bm u(t_{n-1})\big)~.
\end{equation}
First, the focus lies on the last arc $t_{N-1}$ to $t_{N}$. A shift in \ac{TCA} due to a manoeuvre-induced change in the primary's state must be taken into account. This problem has been addressed in prior work (see \cite{armellin2021, pavanello2025c}) by introducing a time variation around the final time step $t_{N}$ with
\begin{equation}
	t_{CA}=t_{N}+\delta t_{CA},
\end{equation}
and extending the final step maps as  
\begin{equation}
\label{eq:PrimaryTcaVar}
	\begin{aligned}	
		\bm x_{p}(t_{CA})&=\mathcal{T}_{\bm x_p}\big(\delta \bm x_{p}(t_{N-1}),~\delta \bm u(t_{N-1}),~\delta t_{CA}\big)\\
		\bm x_{s}(t_{CA})&=\mathcal{T}_{\bm x_s}\big(\delta t_{CA}\big)
	\end{aligned}	
\end{equation}
where the secondary object is only affected by $\delta t_{CA}$ through the natural drift. To handle this time-dependency efficiently within manoeuvre optimisation, an explicit expression of $\delta t_{CA}$ is needed. The definition of closest approach through the orthogonality condition 
\begin{equation}\label{eq:orthCondDA}
	0 \overset{!}{=} \bm r\,(t_{CA}) \cdot \bm v\,(t_{CA})=\mathcal{T}_\textrm{orth}\big(\delta \bm x_{p}(t_{N-1}),\delta \bm u(t_{N-1}),\delta t_{CA}\big)
\end{equation}
allows us to determine $\delta t_{CA}$ for given $\delta \bm x_{p}(t_{N-1}),\delta \bm u(t_{N-1})$. This is done through polynomial inversion facilitated within the \ac{DA} framework used \cite{armellin2021}. The resulting expression for the change in \ac{TCA} due a control action in arc $[t_{N-1}, t_{N}]$ is
\begin{equation}
	\delta t_{CA} = \mathcal{T}_{\delta t_{CA}}\big(\delta \bm x_{p}(t_{N-1}),\delta \bm u(t_{N-1})\big)~.
\end{equation}
Since the motion of both objects in the short-term encounter can be approximated as rectilinear and uniform, truncation to linear maps for the final step is sufficient.
The expression for $\delta t_{CA}$ can be substituted in \cref{eq:PrimaryTcaVar} to express the change in the final primary state as  
\begin{equation}\label{eq:dxpfMap}
	\delta \bm x_{p}(t_{CA})={\mathcal{T}_{\delta \bm x_p}}\big(\delta \bm x_{p}(t_{N-1}),~\delta \bm u(t_{N-1})\big)~,
\end{equation}
effectively incorporating the effect of the \ac{TCA} shift in the polynomial expression of the dynamics. With \cref{eq:dxpMap} and \cref{eq:dxpfMap}, polynomial maps are established to describe the influence of a given control variation on the primary object’s state in each arc. The next step draws inspiration from the early work by Patera and Peterson, where the authors determined the effective thrust direction for a \ac{CAM} by using the gradient of $P_C$ \cite{patera2003}; together with \cref{eq:relStates} and \cref{eq:SMD}, a variation in position and control in each arc $n\in\{0,\dots, N-1\}$ can be mapped to the relative position at \ac{TCA} as
\begin{equation}\label{eq:reachPos}
	\bm r_n(t_{CA})=\mathcal{T}_{\bm r}(\delta \bm x_{p}(t_{n}),~\delta \bm u(t_{n}))~,
\end{equation}
and to the risk metric \ac{SMD} as
\begin{equation}\label{eq:smdMap_n}
	d^2_n=\mathcal{T}_{d^2}(\delta \bm x_{p}(t_{n}),~\delta \bm u(t_{n}))~.
\end{equation}
Note that both quantities are only defined for \ac{TCA}, so $d^2_n$ denotes the updated value due to a manoeuvre in the singular arc with starting node $n$. With the polynomial maps prepared, the \acp{NLP} can be set up.

\subsection{Minimising Lead Time}
The goal is to design \acp{CAM} that are performed as late as possible, thus minimising the manoeuvre lead time. This follows the research direction established in prior work, in which a time-optimal formulation was established and the solution builds on a backward sweep \cite{pavanello2026a}.

A discrete-time \ac{NLP} is set up in which the time horizon \mbox{$t\in\left[-t_W,0\right]$} is split into $N$ arcs. The lead time $t_L$, as the start of the first arc with non-zero control, is used as the cost function. Polynomial maps capture the dynamics relative to the nominal trajectory at each node $n$ as
\begin{equation}\label{eq:NlpDynamics}
\delta \bm{x}_{p}(t_{n})  =  \mathcal{T}_{\delta x_p}\big(\delta \bm{x}_{p}(t_{n-1}), \mathcal{C}_{f,\sigma}(\bm{p}_{n-1})\big)~,
\end{equation}
with $\mathcal{C}_{f,\sigma}$ being the control mapping defined in \Cref{eq:pseudoToAccel}.
Therefore, the free lead time formulation of the \ac{NLP} is
\begin{subequations}
\begin{align}
    \min\limits_{\bm{p}_n} \quad & t_L \\
    \mbox{s.t.} 
    \quad &  \text{\cref{eq:NlpDynamics}}   \\
    \quad & \delta\bm{x}_{p}(-t_L) = \bm{0} \label{eq:nlpLead:init}\\
    \quad & \bm{p}_\textrm{min} \le \bm{p}_n \le \bm{p}_\textrm{max} \label{eq:nlpLead:ctr}\\
    \quad & {d}^2 \ge {d_{\mathrm{min}}}^2, \label{eq:nlpLead:smd}
\end{align}
\label{eq:nlpLead}
\end{subequations}
with zero initial deviation at the start of the first control arc \cref{eq:nlpLead:init}, constraints on the control \cref{eq:nlpLead:ctr} and the risk threshold at \ac{TCA} captured through \ac{SMD} via \cref{eq:nlpLead:smd}. 
By solving, we find the lead time $t_L$ as the time prior to \ac{TCA} at which the first control of the \ac{CAM} needs to be performed to still safely avoid the encounter (cf. \cref{fig:ConjScenario}). 

\subsubsection{The Minimum-Lead Time Method}
The developed solution method to solve \mbox{Problem \ref{eq:nlpLead}} is named \ac{MTCOLA}. The core idea of \ac{MTCOLA} is built on a backward sweep approach: starting from \ac{TCA}, the control profile is updated by stacking actions backward in time until the risk threshold is reached. This is achieved by using polynomial maps of the state and cost functions to derive locally optimal control solutions in each arc.

We begin at the last arc, where a variation in position and control at $t_{N-1}$ is mapped to the \ac{SMD} of the last time step $t_{CA}$ as 
\begin{equation}\label{eq:smdMapNm1}
	d^2_{N-1}=\mathcal{T}_{d^2}(\delta \bm x_{p}(t_{N-1}),~\delta \bm u(t_{N-1}))~.
\end{equation}
On the basis of $d_{N-1}^2$ and a given change in position $\delta \bm x_{p}(t_{N-1})$, we can choose a control $\delta \bm u^{*}$ that is constrained using the polynomial expression in \cref{eq:constraint}. The constrained control is used to compute the change in relative position $\bm r\,(t_{CA})$ by evaluating \cref{eq:smdMap_n}. Then, an updated $P_C$ for the adjusted relative position is calculated using Chan's formula \cite{chan2018}. If the safe threshold in $P_{C}$ is reached, the \ac{CAM} is deemed effective, and the iteration is stopped. Otherwise the maps are prepared for the next step by linking the initial position $\delta \bm x_{p}(t_{N-2})$ and control $\delta \bm u(t_{N-2})$ of the next iteration's arc (backwards in time) to the position of the current step $\delta \bm x_{p}(t_{N-1})$. From substituting with 
\begin{equation}
	\delta \bm x_{p}(t_{N-1}) = \mathcal{T}_{\delta \bm x_{p}}(\delta \bm x_{p}(t_{N-2}),~ \delta \bm u(t_{N-2}))~,
\end{equation}
the next arc's state is 
\begin{equation}
	\delta \bm x_{p}(t_{CA}) = \mathcal{T}_{\delta \bm x_{p}}(\delta \bm x_p(t_{N-2}),~\delta \bm u(t_{N-2}),~\delta \bm u^{*}(t_{N-1}))~,
\end{equation}
and the \ac{SMD} formulations
\begin{equation}
	d^2_{N-2} = \mathcal{T}_{d^2}(\delta \bm x_p(t_{N-2}),~\delta \bm u(t_{N-2}),~\delta \bm u^{*}(t_{N-1}))~,
\end{equation}
and can be used to find the next optimal control $\delta \bm u^{*}(t_{N-2})$. Crucially, both maps differ from those defined in \cref{sec:polynomialMaps}, in that they inherit the effect of the applied control $\delta \bm u^{*}(t_{N-1})$. Following the approach in \cite{pavanello2026a}, iterations are performed in a backward sweep until the risk threshold is surpassed. The obtained solution is only time-optimal with regard to the chosen time discretisation. Although a refinement of the last arc has been proposed, this feature is not included in this work. Therefore, the manoeuvre solutions overshoot the risk threshold based on the choice of $N$.

In the original formulation, the selected control action at each time step is the one that maximizes the metric rate of change at \ac{TCA}. However, with the nominal configuration as minimum drag, no control deviations in positive $\bm{\hat x}_W$-direction can be applied. As such, strategies that greedily apply the control corresponding to the risk gradient are not applicable, since movement on the B-plane is constrained. We will therefore focus next on deriving an approximation for the B-plane positions reachable via a constrained manoeuvre to analyse possible strategies of aerodynamic \acp{CAM} and find a better rationale for choosing suitable controls.

\subsubsection{Reachable Set Boundary Computation}
The goal of this section is to develop a method to support the selection of optimal controls in the \ac{MTCOLA} approach. The idea is to relate available aerodynamic accelerations in three dimensions to their effect on the conjunction. We achieve this by identifying the maximum reachable position change in various directions on the B-plane (two-dimensional projection) and analysing how these evolve with respect to the control constraints and the manoeuvre timing. 

\begin{figure}[tb]
	\centering
	\includegraphics[clip, trim={0 2cm 0 0},width=1\figwidth]{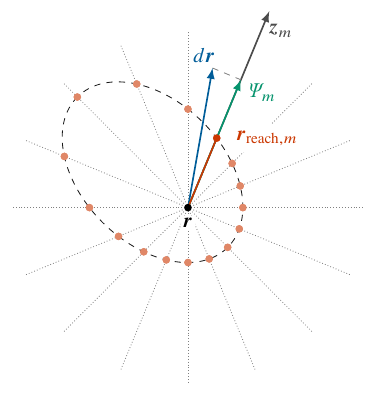}
	\caption{Illustration of reachable boundary approximation on the B-plane.}
	\label{fig:reachSet}
\end{figure}

To analyse reachable positions on the B-plane induced by aerodynamic control, a simplified local method is employed. Instead of solving an optimal control problem to construct the full reachable set over the entire time horizon as commonly done in the literature \cite{bansal2017}, this approach directly approximates the reachable set boundary over a single discrete time step. The relative B-plane position $\bm r_n(t_{CA})$ at closest approach for a deviation in state and control at $t_n$ is the polynomial map in \cref{eq:smdMap_n}. Next, $M$ rays $\bm z_{m}\in\mathbb{R}^{2}$ (in B-plane coordinates) starting at the current position and with equal angular distance are constructed as
\begin{equation}
	\bm z_{m} = \begin{bmatrix} \cos(\theta_m) \\ \sin(\theta_m) \end{bmatrix}, \quad \theta_m = \frac{2\pi \cdot m}{M}, \quad m = 0, \ldots, M-1~.
\end{equation}
The goal is to determine the control deviation $\delta \bm u_{m}^{*}$ that maximises the change of the relative position in the direction of each ray. Since in the following we consider the control effect of a single arc, the time reference $t_n$ is omitted for clarity. Using the Jacobian of the relative position polynomial map with respect to the control
\begin{equation}
	\bm J(\bm u) := \frac{\partial \bm r}{\partial \bm u} \in \mathbb{R}^{2\times 3}~,
\end{equation}
the change in relative position $d \bm r$ for a small control increment $d \bm u$ can be expressed with the first-order approximation
\begin{equation}
	d \bm r = \bm J(\bm u) \, d \bm u~.
\end{equation}
The projection $\varPsi_{m}$ onto the ray $\bm z_{m}$ is the dot product as
\begin{equation}
	\varPsi_{m}(d \bm u) = (d \bm r)^{T}\bm z_{m} \approx \big(\bm J(\bm u) ~ d \bm u\big)^{T}\bm z_{m}~.
\end{equation}
The instantaneous maximisation of $\varPsi_{m}(d\bm u)$ yields the desired control direction
\begin{equation}
	\delta \bm u^{*}_{m} = \frac{\bm J(\bm u)^T \bm z_{m}}{\|\bm J(\bm u)^T \bm z_{m}\|}~.
\end{equation}
Note that $|d \bm u| = 1$ is assumed since only the control direction is required. \cref{fig:reachSet} shows the construction of the projection $\varPsi_{m}$ onto one ray using an arbitrary $d \bm r$ for illustration purposes. The constrained control magnitude for the desired direction $\delta \bm u^{*}_{m}$ is obtained using the aerodynamic constraint in \cref{eq:constraint} and used to evaluate the reachable relative position $\bm r_{\mathrm{reach},m}$ from \cref{eq:reachPos}. The resulting set of $M$ controls and the corresponding change in relative position in each of the $N$ arcs of the \ac{CAM} is used to get an intuition of the dynamics on the B-plane and their dependence on the control authority. Since sampling occurs over the one-dimensional Euclidean angle and pre-computed polynomial maps are used, this approach provides an efficient approximation of the reachable boundary. The method can be applied to any underlying dynamics, as the polynomial maps are computed via \ac{DA}, and it can also be applied to propulsive \acp{CAM}.

\subsubsection{Control Authority of Aerodynamic \acp{CAM}}
\begin{figure}[tb]
	\centering
	\begin{subfigure}[b]{1\figwidth}
		\centering
		\includegraphics[width=\textwidth]{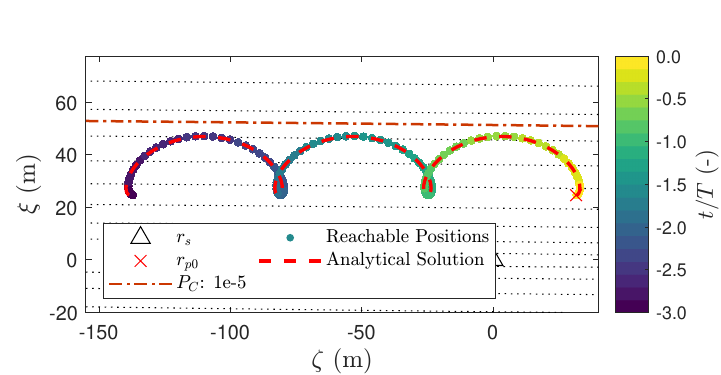}
		\caption{Drag-only, specular surface, analytical solution from \cite{hernando-ayuso2021}}	
		\label{fig:reachDrag}
	\end{subfigure}  
	\begin{subfigure}[b]{1\figwidth}
		\centering
		\includegraphics[width=\textwidth]{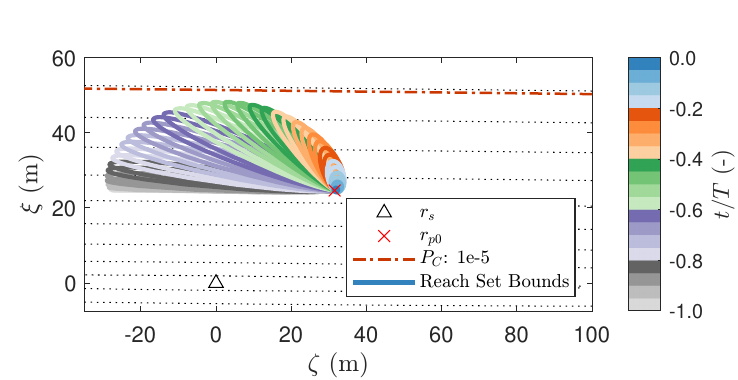}
		\caption{Lift-based, specular surface}
		\label{fig:reachLift}
	\end{subfigure}	
	\begin{subfigure}[b]{1\figwidth}
		\centering
		\includegraphics[width=\textwidth]{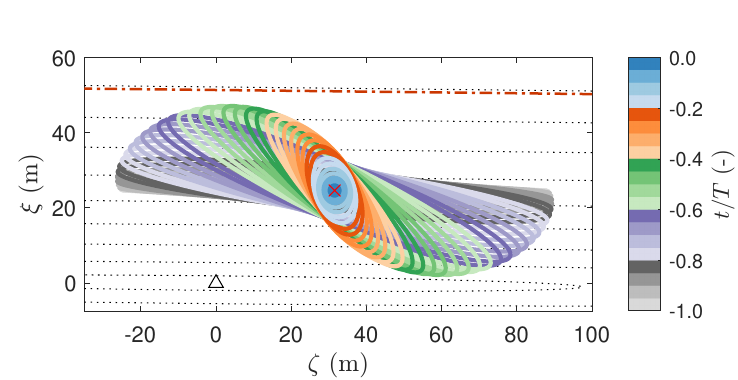}
		\caption{Propulsive}
		\label{fig:reachProp}
	\end{subfigure}
	\caption{Reachable set boundaries for one-arc manoeuvres with length $T/40$ and varying timing $t/T$ for the three control constraints. The encounter is at \qty{300}{\kilo\metre} and $\kappa=115$ with iso-probability ellipses for $I_{CDM}=77$ in \cref{tab:mtColaScenarios}. Note that (a) shows 120 arcs over three orbits, while (b) and (c) show 40 arcs over one orbit.}
	\label{fig:reachability}	
\end{figure}

To develop an understanding of the issue and inform the development of algorithms for lift-based actuation compared to drag-only and propulsive \acp{CAM}, we analyse the reachable positions on the B-plane for varying control authorities, manoeuvre timings and $\sigma$ values. The reachable boundaries of drag-only, lift-based, and propulsive manoeuvres are displayed in \cref{fig:reachability}. They show the B-plane positions reachable via one single-arc manouevre. The timing of the manoeuvre is given in orbital revolutions prior to \ac{TCA} by relating it to the orbital period $T$ of the primary object. Both objects are in quasi-circular orbits with mutual inclination $\kappa=\ang{115}$ and an altitude of \qty{300}{\kilo\metre}. \cref{fig:reachDrag} shows reachable positions via in-track accelerations due to a drag-only manoeuvre. Each dot represents the position after a singular arc at time $t/T$ with maximum drag, followed by drift in the nominal drag configuration until \ac{TCA}. The timing of the manoeuvre is varied over $120$ arcs in three orbits prior to closest approach. Results of the in-track \acp{CAM} align with the manoeuvre effects presented by Aida \cite{aida2016} and replicate the analytical expression for tangential B-plane kinematics derived by Hernando-Ayuso and Bombardelli \cite{hernando-ayuso2021}, shown as the dashed line in \cref{fig:reachDrag}. 

Several key observations emerge. Arcs further from \ac{TCA} tend to produce larger absolute position changes, but effectiveness is not monotonic. Some arcs just under one full orbital revolution exhibit larger position changes than those slightly over one revolution. The drag manoeuvres also reveal a distinctive characteristic of the B-plane: changes to the in-track separation distance, achieved through tangential manoeuvres close to $t/T=1$, cause a deviation in $\zeta$, thus also called the timing axis. The radial separation of both objects at \ac{TCA} is reflected as a shift in $\xi$, which for tangential manoeuvres is most effective close to $t/T=0.5$. Along with the common elongation of the iso-probability ellipses along $\zeta$, this creates a critical timing effect: performing the same tangential manoeuvre earlier does not necessarily result in a greater reduction in risk. For example, by superposition of individual discrete manoeuvres, a two-arc manoeuvre near $t/T=0.5$ will cross the $P_C$ threshold, while the same manoeuvre near $t/T=1$ will not. \ac{CAM}s close to half a revolution therefore align better with the collision probability gradient of this conjunction, that is, perpendicular to the dotted iso-probability ellipses. This again demonstrates effectively why along-track separation is not a suitable metric for efficient \ac{CAM} design since it fails to capture the true nature of collision risk. The described relationships are already established in literature for tangential \acp{CAM}. However, the reachable set boundary computation proves especially useful in increasing intuition of complex constraints and a variation in maneouvre timing, and it could even be applied to non-circular orbits.

\cref{fig:reachLift} demonstrates the impact of lift-based manoeuvres on the reachable boundary over 40 arcs spanning one orbit. The number of rays was set to $M=500$ for visualisation. Each oval of a colour represents the reachable boundary for a manoeuvre performed over one arc, followed by uncontrolled drift. Aerodynamic lift expands the reachable boundary from a line (drag-only) to an oval, whose width scales with the available $L/D$ ratio. The propulsive case in \cref{fig:reachProp} helps to illustrate the directional constraint of drag-dominated manoeuvres: Aerodynamic \acp{CAM} starting from the minimum-drag configuration as set up in this work may only increase drag, thus decreasing the \ac{SMA} compared to the nominal trajectory. Therefore, movement on the B-plane is restricted in the positive $\xi$ direction. This is not the case for the propulsive \ac{CAM} since it is able to increase \ac{SMA}. The same goes for the timing axis $\zeta$: increased drag leads to the encounter being earlier (negative change in $\zeta$), and more so if performed longer prior to \ac{TCA}, due to the decrease in \ac{SMA} reducing the orbital period of the primary. Tangential accelerations essentially create phasing effects that elongate the lift-based reachability ovals along the $\zeta$-axis for manoeuvres applied more distant from \ac{TCA}. The direction of maximum deviation in the lift-based case remains within the top-left quadrant for all lead times $t_L>0.2$, and the cyclic variations in direction decrease for longer lead times. These directional relationships assume the radially inward $\xi$-axis orientation, used throughout this work.

\begin{figure}[tb]
	\centering
	\includegraphics[width=\figwidth]{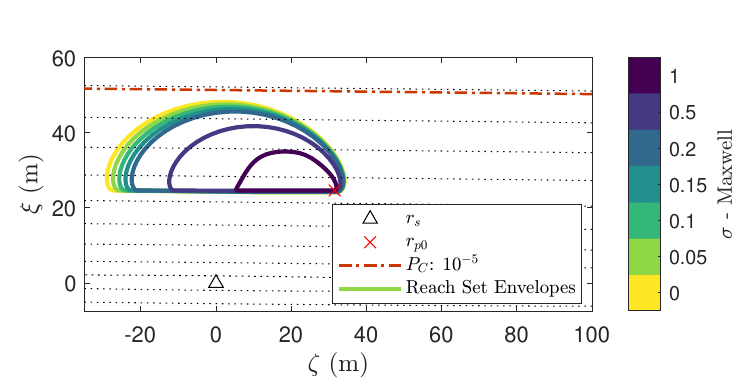}
	\caption{Reachable set envelopes for one-arc lift-based manoeuvres with length $T/40$ at any point of one orbit for varying surface properties $\sigma$.}
	\label{fig:reachHullLift_sigma}	
\end{figure}

\begin{figure}[tb]
	\centering
	\includegraphics[width=0.5\textwidth]{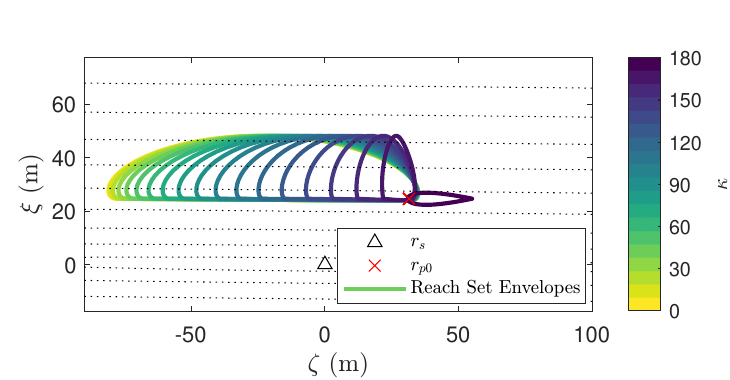}
	\caption{Reachable set envelopes for one-arc lift-based manoeuvres with length $T/40$ at any point of one orbit for varying mutual inclinations $\kappa$.}
	\label{fig:reachHullLift_kappa}
\end{figure}

\begin{figure}[tb]
	\centering
	\includegraphics[width=\figwidth]{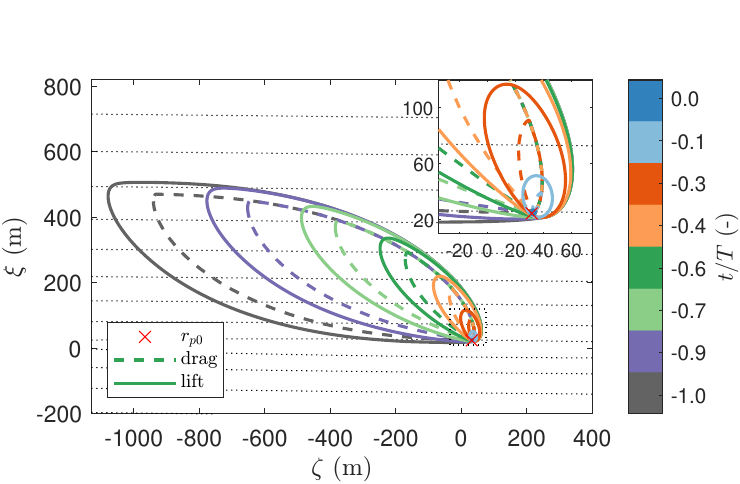}    
	\caption{Reachable set hulls for multi-arc manoeuvre profiles and varying lead times. Comparison of lift and drag under specular reflection ($\sigma=1$).}
	\label{fig:reachHull}	
\end{figure}

Next, we construct the full set of maximum positions reachable via one-arc lift-based manoeuvres at any time during an orbit prior to \ac{TCA}. This is essentially the envelope of all the reachable set boundaries shown in \cref{fig:reachLift}. The influence of $\sigma$ on these envelopes is displayed in \cref{fig:reachHullLift_sigma}. Clearly, near-specular properties inhibit the largest possible position changes. \Cref{fig:reachHullLift_kappa} shows the effect of the mutual inclination on the envelopes. For low $\kappa$ the reachable positions on the B-plane exhibit the largest $\zeta$-component, while for increasing $\kappa$ a change on the timing axis becomes harder to achieve. The head-on case ($\kappa=\ang{180}$) again stands out since the timing has no effect on the encounter, giving the hull an entirely different shape.
Since the reachable sets at each manoeuvre timing are convex, the Minkowski sum of these sets over one orbit yields the convex hull of the complete reachable set for the entire manoeuvre duration. This allows the reachable set for a complex control trajectory to be approximated for the given discretisation and compared between drag and lift, see \cref{fig:reachHull}. Here, convex hulls for varying warning times below one orbit are contrasted for drag (dashed) and lift (solid) manoeuvres under specular reflection. Lift-based manoeuvres achieve larger position changes in all directions except $-\xi$, demonstrating their clear advantage over drag-only approaches.
While no prior work known to the authors has visualised aerodynamic lift effects in this manner, Bombardelli and Hernando-Ayuso demonstrated reachability analysis for impulsive manoeuvres based on the linearised B-plane kinematics \cite{bombardelli2015}. Their resulting reachable domains take elliptical shapes that vary with manoeuvre timing and resemble the effects described and seen in \cref{fig:reachProp}. The reachable set computation method developed for this work makes no assumption of circular orbits and can be extended to approximate long manoeuvres with varying control directions by superposition of individual reachable boundaries. The shown visualisation helps clarify the relation of control authority, conjunction geometry, and collision risk reduction that is not exclusive to aerodynamic manoeuvres and low lead times. Next, the reachable sets are used to find the best control in minimum lead time scenarios.

\subsubsection{Heuristic Control Strategy}
\begin{figure}[tb]
	\centering
	\includegraphics[clip, trim={0 0 0 0},width=0.9\figwidth]{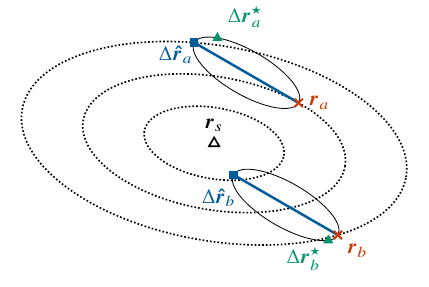}	
	\caption{Illustration of the reachable set control strategy on the B-plane. Shown are three levels of \ac{SMD} as dotted lines, primary starting positions $\bm r_{a}$ and $\bm r_{b}$, the reachable set boundary for control in one arc, and primary after the manoeuvre of maximum absolute deviation $\Delta\bm{\hat r}$ and maximised \ac{SMD} change $\Delta\bm{r}^\star$.}
	\label{fig:reachControl}
\end{figure}

\begin{figure}[tb]
	\centering
	\includegraphics[clip, trim={0 0 0 0},width=0.9\figwidth]{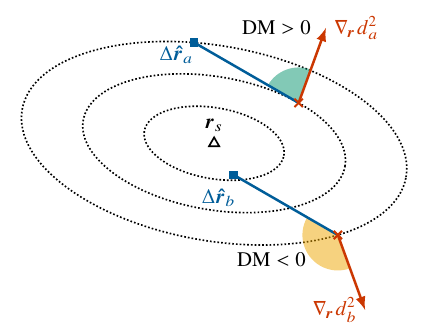}		
	\caption{Illustration of descent metric construction on the B-plane with three levels of \ac{SMD} as dotted lines, the dashed switching line, the gradients in two points $\nabla_{\bm r} d^2$, the direction of maximum absolute deviation $\Delta\bm{\hat r}_{CA}$ and the risk descent metric $\textrm{DM}$.}
	\label{fig:descentMetric}
\end{figure}

The reachable set boundary expresses information about the maximum position change on the B-plane and is used during runtime within \ac{MTCOLA}. Instead of choosing the control gradient direction and accepting that the aerodynamic constraint might limit the manoeuvre effect to a minimum, we include the constraint via the reachable sets. For this, we find the control corresponding to the highest reachable \ac{SMD} in each control arc out of all $M$ rays, with
\begin{equation}
	d^{2\,\star}= \max_{m} d^2_m~,
\end{equation}
and the control that corresponds to the maximum change in position on the B-plane, which is
\begin{equation}
	\Delta{\hat r} = \max_{m}~~ \lVert \bm r_{\mathrm{reach},m}- \bm r \rVert~.
\end{equation}
In \cref{fig:reachControl}, two examples of reachable boundaries over one control arc are shown, starting at different B-plane positions. For purposes of illustration, the starting positions of the primary $\bm{r}_a,~\bm{r}_b$ are chosen to be on iso-probabilty ellipses. The reachable positions corresponding to a maximum absolute distance are $\Delta\bm{\hat r}_a,~\Delta\bm{\hat r}_b$, while $\Delta\bm{r}_a^\star,~\Delta\bm{r}_b^\star$ are the positions that maximise \ac{SMD}. Clearly, in case (a) the point corresponding to the largest \ac{SMD} increase $\Delta\bm{r}_a^\star$ is the best option. On the other hand, in (b) this is no longer the case. Although choosing the control for $\Delta\bm{r}_b^\star$ leads to a slight increase in \ac{SMD}, this gain is minimal. While the manoeuvre for $\Delta\bm{\hat r}_b$ decreases \ac{SMD} (increasing collision risk) in that arc, it enables faster risk reduction in future control arcs, resulting in an overall shorter manoeuvre lead time. In such cases, it is beneficial to apply a heuristic control method that permits the \ac{SMD} to decrease for some iterations until a switching line is crossed and the control corresponding to maximum risk reduction is effective again.  This approach relies on the reachable boundaries staying similar in shape and orientation over the following control arcs, as was shown in \cref{fig:reachability}. Two properties can be identified:
\begin{enumerate}
	\item the direction towards $\Delta\bm{\hat r}$ varies over the orbit position relative to \ac{TCA} but remains within one B-plane quadrant for reasonable lead times, and
	\item the achievable B-plane displacements are small but non-zero near \ac{TCA} and increase for arcs executed earlier.
\end{enumerate}
While no formal proofs are provided, these properties were observed in all analyses performed. Leveraging these properties, we introduce a risk descent metric ($\textrm{DM}$) that evaluates the alignment between the gradient direction of the \ac{SMD} in the B-plane, $\nabla_{\bm r} d^2$, and the direction of maximum achievable position change from the reachable set analysis, $\Delta\hat{\bm r}$, with
\begin{equation}
	\textrm{DM} = (\nabla_{\bm r} d^{2})^{T} \Delta\hat{\bm r}~.
\end{equation}
When $\textrm{DM} \ge 0$, the aerodynamic control constraint allows movement in a favourable direction that leads to an increase in \ac{SMD}, and a decrease in $P_C$. When $\textrm{DM} < 0$, the control constraint renders ascent in \ac{SMD} ineffective and the control corresponding to the maximum absolute position change $\bm \hat r$ is applied instead. The heuristic control strategy is
\begin{equation}
	DM\quad
	\begin{cases}
		\ge 0~~\text{choose control for}~\Delta\bm{r^\star} \\
		< 0~~\text{choose control for}~\Delta\bm{\hat r}~.
	\end{cases}
\end{equation}
Both cases are illustrated in \cref{fig:descentMetric}. The switching line of $\textrm{DM} = 0$, which depends on the $\Delta\bm{\hat r}$ direction, is displayed as a dashed line. This approach essentially permits strategic detour moves that will, for some iterations, increase collision probability but ultimately enable more effective manoeuvres. 

Unlike greedy gradient-based control, this approach naturally accounts for the directional constraint inherent to lift-based \acp{CAM}: since $\Delta\bm{\hat r}$ remains within a single B-plane quadrant throughout the encounter (top-left in this scenario), with cyclic variations that narrow for longer lead times (cf. \cref{fig:reachability}), the temporary 'detours' occurring when $\textrm{DM} < 0$ remain bounded in direction. Convergence follows because every applied control produces a non-zero, directionally bounded displacement that advances toward the iso-probability ellipse of the required $P_C$, which is therefore reached after a finite number of control arcs.

\begin{figure}[t]
	\centering
	\includegraphics[clip, trim={0 0 0 0},width=1\figwidth]{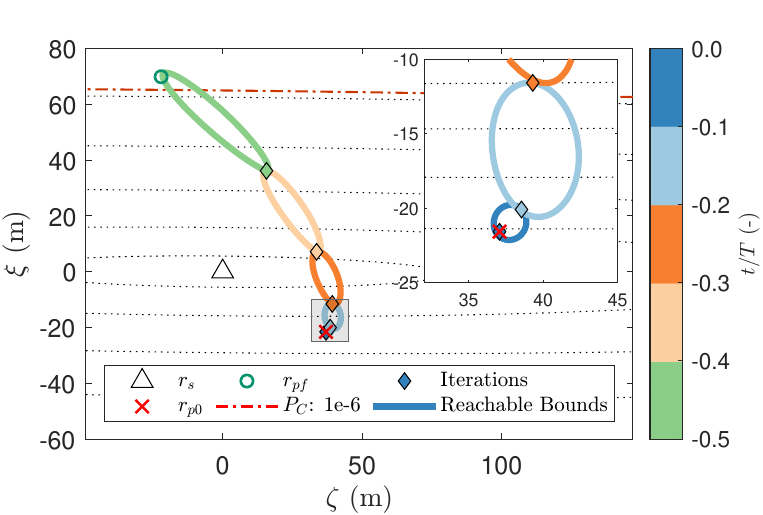}
	\caption{Visualisation of the position evolution and the reachable set boundary throughout a \ac{CAM} within \ac{MTCOLA} for 5 arcs.}
	\label{fig:reachVis}
\end{figure}

\cref{fig:reachVis} demonstrates how the heuristic approach is applied in iterations, stepping back in time from \ac{TCA}, selecting the optimal control direction in each arc. In every iteration, the point on the reachable boundary of maximum \ac{SMD} is chosen (here $DM>0$), and the corresponding control is selected. For the next iteration, this point becomes the anchor from which the next reachable boundary is computed. In this work, the number of rays used for the reachable boundary computation is set to $M=100$, to achieve a balance between computational cost and discretisation fidelity. \ac{MTCOLA} proves robust across all the analysed scenarios. Due to its deterministic nature, it always finds a solution for the given discretisation, however, optimality is not guaranteed. 

\subsection{Minimising Decay}
Aerodynamic \acp{CAM} as set up in this work inevitably apply drag forces that induce increased orbital decay not accounted for in \ac{MTCOLA}. Minimising this effect, however, can be solved with more common direct single shooting methods that use the polynomial maps of \cref{sec:polynomialMaps}. The manoeuvre lead time is not optimised for and may reach up to the warning time with $t_L=t_W$. Instead, we aim to minimise the decay of the primary spacecraft due to the \ac{CAM} as 
\begin{equation}\label{eq:dSMA}
	\Delta a = a_{p,\textrm{CAM}}(t_{CA})-a_{p,\textrm{NOM}}(t_{CA})~,
\end{equation}
with the \ac{SMA} of the primary spacecraft before $a_{p,\textrm{NOM}}$ and after the manoeuvre $a_{p,\textrm{CAM}}$. The \ac{NLP} is

\begin{subequations}
\begin{align}
    \min\limits_{\bm{p}_n} \quad & -\Delta a = f\left(\delta \bm{x}_p(t_{CA})\right)  \\
    \mbox{s.t.} 
    \quad &  \text{\cref{eq:NlpDynamics}}   \\
    \quad & \bm{p}_\textrm{min} \le \bm{p}_n \le \bm{p}_\textrm{max} \\
    \quad & {d}^2 \ge {d_{\mathrm{min}}}^2~,
\end{align}
\label{eq:NLP_minDecay}
\end{subequations}
with constraints on the control in \cref{eq:nlpLead:ctr} and \ac{SMD} in \cref{eq:nlpLead:smd}. We solve this using the direct single shooting approach. Polynomial maps for the primary position are created by stepping forward in time for all arcs starting at $t_W$. The primary's deviation from the nominal trajectory is fully described via \cref{eq:dxpMap} for $n\in\{0,\dots,N-2\}$ and \cref{eq:dxpfMap} for the arc $n=N-1$ that inherently tracks \ac{TCA}.

Using the fundamental equation of the specific energy of the Keplerian orbit, the decay in \ac{SMA} becomes
\begin{equation}\label{eq:dSMA2}
	\Delta a = 	\left(\frac{2}{r_{p,\textrm{CAM}}} - \frac{v_{p,\textrm{CAM}}^2}{\mu_\textrm{E}}\right)^{-1}
	- \left(\frac{2}{r_{p,\textrm{NOM}}} - \frac{v_{p,\textrm{NOM}}^2}{\mu_\textrm{E}}\right)^{-1}~,
\end{equation}
with 
\begin{equation}
	\begin{aligned}
		r_\textrm{p,CAM}&=\lVert \bm r_{p,\textrm{NOM}}(t_{CA})+\delta \bm r_p(t_{CA})\rVert \\
		v_{p,\textrm{CAM}}&=\lVert\bm v_{p,\textrm{NOM}}(t_{CA})+\delta \bm v_p(t_{CA})\rVert~,
	\end{aligned}
\end{equation}
and the nominal position $\bm r_{p,\textrm{NOM}}(t_{CA})$ and velocity $\bm v_{p,\textrm{NOM}}(t_{CA})$ at \ac{TCA} without a \ac{CAM} performed. The resulting cost function is expressed as a polynomial map. Importantly, the maps are interlinked through the state deviations $\delta \bm x_p$. This means that the effect of the first arc control on the final state and thus decay and collision risk (\cref{eq:smdMapNm1}) is captured through the state transition constraint in the \ac{NLP}. By using the chain rule in the \ac{DA} polynomial, given that state deviations are only dependent on control actions that are applied at the previous time node, it is possible to write the \ac{SMA} decay as a function of the sequence of control actions
\begin{equation}
	\Delta a 
	= 
	\mathcal{T}_{\Delta a}
	\big(
	\mathcal{C}_f(\bm p_0),
	\mathcal{C}_f(\bm p_1),
	\dots,
	\mathcal{C}_{f}(\bm p_{N-1})
	\big)~,
	\label{eq:dSMA_map}
\end{equation}
where $\mathcal{C}_f(\bm p_n)$, as defined in \cref{eq:pseudoToAccel}, maps the pseudo-control in arc $n$ to its corresponding aerodynamic acceleration. By this, $\Delta a$ is an explicitly differentiable function of all $\bm p_n$. The polynomial map representation also ensures tracking of \ac{TCA} within the last time step, avoiding the need for an additional nonlinear constraint. Furthermore, the maps ensure comparability with the \ac{MTCOLA} approach, which can be used to provide a reliable initial guess for a minimum lead time solution that is optimised for decay via direct shooting.

% -------------------------------------------------------------------------------------%
% -------------------------------------------------------------------------------------%
% ------------------------------------- Results --------------------------------------- %
% -------------------------------------------------------------------------------------%
% -------------------------------------------------------------------------------------%

\section{Aerodynamic \ac{CAM} Performance}\label{sec:results}
In the following, the simulation setup and \ac{CDM} data used is briefly described followed by results and comparisons of lift-based, drag-only and propulsive manoevres.
%The code for the presented methods and the data from which the results are drawn are available online and can be found at \cite{??} \textcolor{red}{todo, make repo and data available}. 
The setup of \ac{DA} maps for the direct shooting approach and the \ac{MTCOLA} algorithm is implemented as C\texttt{++} routines based on the \textit{DACE} framework \cite{rasotto2016,wittig2025}. The creation of scenario configurations, implementation of the single shooting method, and post-processing is done within MATLAB. All data is stored in \textit{HDF5} format to ensure modularity and structured storage of simulation runs. The direct shooting approach is implemented using \textit{CasADi} in MATLAB \cite{andersson2019} and solved using the \textit{Interior Point OPTimizer} (\textit{IPOPT}) \cite{ipopt} (version 3.14.11) with the linear solver \textit{MUltifrontal Massively Parallel sparse direct Solver} (\textit{MUMPS}) \cite{mumps} (version 5.4.1). The polynomial structure of the \ac{NLP} provides all derivatives up to order $k=2$ and eliminates the need for automatic differentiation within \textit{CasADi}.
All simulations were performed on an AMD Ryzen 7 7730U processor using a single core. The MATLAB components ran on Windows 11 with \qty{16}{\giga\byte} RAM (MATLAB R2023b), while the C\texttt{++} modules were executed within Windows Subsystem for Linux (WSL) with 8GB RAM allocated. 

The \ac{MTCOLA} method demonstrates reliable convergence and moderate computational efficiency. It finds a minimum lead time \ac{CAM} with $N=200$ arcs spanning two orbital periods in less than \qty{0.5}{\second} including reading and writing of in- and output data. The direct shooting method typically converges in less than \qty{3}{\second} for scenarios with $N \le 50$ arcs. However, the setup of extracting and formulating the polynomial maps as \textit{CasADi} symbolic expressions requires \qty{3}{\second} as well, doubling the total time to run a single scenario. Computational times of competitive numerical methods  for \ac{CAM} design, such as sequential convex programs or polynomial programs, are shown to be moderately faster \cite{pavanello2026a,pavanello2025c}, while (semi-)analytical methods are much faster \cite{hernando-ayuso2021,devittori2022}. However, the current implementation has not been optimised in terms of computational speed. It still enabled thousands of optimisations to be performed for parametric analyses of lift-based \acp{CAM} at varying altitudes, warning times, covariance magnitudes and encounter geometries. The following analyses are based on selected conjunction scenarios derived from the \textit{ESA Collision Avoidance Challenge} \cite{esa2019} dataset listed in \cref{tab:mtColaScenarios}. Since the original data does not contain sufficient conjunctions at very low altitudes, orbit heights were modified to emulate \acs{VLEO} altitudes between \qty{200}{\kilo\metre} and \qty{400}{\kilo\metre} depending on the analysis.

\subsection{Minimum Lead Time Performance}
\begin{figure}[t]
	\centering
	\begin{subfigure}[b]{1\figwidth}
		\centering
		\includegraphics[width=0.99\textwidth]{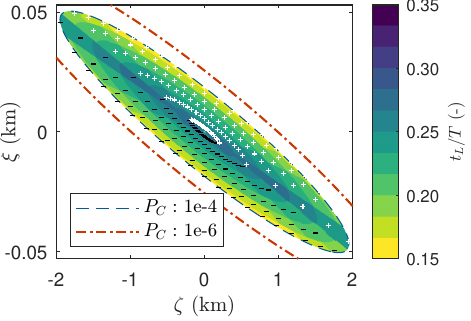}
		\caption{Propulsive}
		\label{fig:perfProp}
	\end{subfigure}
	\begin{subfigure}[b]{1\figwidth}
		\centering
		\includegraphics[width=0.99\textwidth]{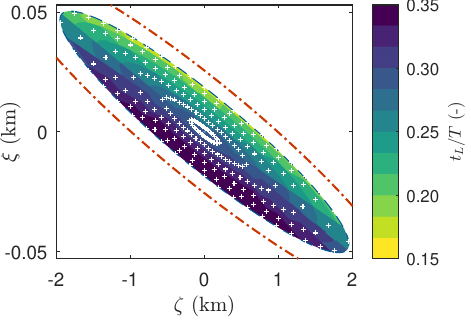}
		\caption{Lift-based}		
		\label{fig:perfLift}        
	\end{subfigure}
	\caption{Performance comparison of propulsive and lift-based minimum lead times for \ac{CAM}s from \ac{MTCOLA}. The altitude is \qty{350}{\kilo\metre}, $\kappa=\ang{160}$, $\sigma=0$ and $I_\textrm{CDM}=1$. For each case, 240 starting positions were analysed and the white plus and black minus signs at the starting positions indicate the direction of travel in positive and negative $\xi$-direction to reach the $P_C$ threshold.}
	\label{fig:performanceMinLead}        
\end{figure}

Minimum lead times for drag-only, lift-based, and propulsive \acp{CAM} were assessed using \ac{MTCOLA}. To highlight the directional nature of aerodynamic \acp{CAM} as set up in this study, the performance of propulsive and lift-based manoeuvres is shown over a grid of 240 initial B-plane positions in \cref{fig:performanceMinLead}. The direction of $\xi$-axis travel to reach the safe threshold is marked via the black minus and white plus signs at each starting position to highlight the directional constraint of aerodynamic \acp{CAM}. For the propulsive baseline in \cref{fig:perfProp}, the achieved lead times are longest near the major axis of the dashed iso-probability ellipse, as the required distance to reach the risk threshold is greater. Lift-based \acp{CAM} (see \cref{fig:perfLift}) exhibit directional limitations: since all the trajectories lead towards the positive $\xi$-direction, the lower part of the initial position sampling ellipse has higher lead times.

\begin{figure}[tb]
	\centering
	\begin{subfigure}[t]{1\columnwidth}
		\centering
		\includegraphics[height=4cm]{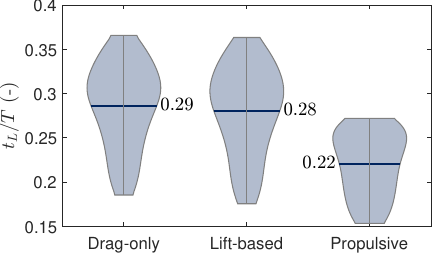}
		\vspace{1mm}
		\caption{Lead time}		
		\label{fig:performanceComparison_position_leadTime}
	\end{subfigure}
	\begin{subfigure}[t]{1\columnwidth}
		\centering
		\includegraphics[height=4cm]{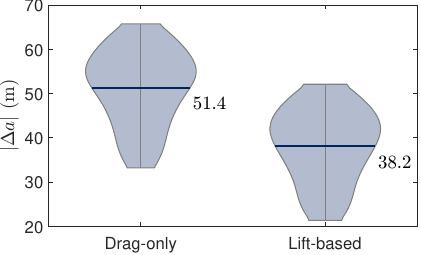}
		\vspace{1mm}
		\caption{Decay in semi-major axis}		
		\label{fig:performanceComparison_position_decay}
	\end{subfigure}
	\caption{Performance distribution of drag-only, lift-based and propulsive \acp{CAM} from the \ac{MTCOLA} method and for 240 B-plane starting positions.  Altitude is \qty{350}{\kilo\metre}, $\kappa=\ang{160}$, $\sigma=0$ and $I_\textrm{CDM}=1$.}
	\label{fig:performanceComparison_position}
\end{figure}

A direct comparison of drag-only, lift-based, and propulsive minimum lead time \acp{CAM} for the 240 analysed positions is shown in \cref{fig:performanceComparison_position}.  The violin plots indicate the statistics of the analysed data, with the width along the $x$-axis representing the probability density at the corresponding $y$-value. Mean values are represented by horizontal lines. While the propulsive approach outperforms the aerodynamic methods in terms of lead time, the lift-based method achieves slightly lower lead times and significantly reduced decay compared to the drag-only baseline. Although \ac{MTCOLA} does not specifically aim to limit orbital decay, a reduction in decay is observed for lift-based \acp{CAM}. This is because control using lift components achieves greater positional changes on the B-plane (cf. \cref{fig:reachHull}), resulting in lower required drag forces and therefore reduced decay.

\begin{figure*}[tb]
	\centering
	\begin{subfigure}[t]{0.7\textwidth}
		\centering
		\includegraphics[width=0.99\textwidth]{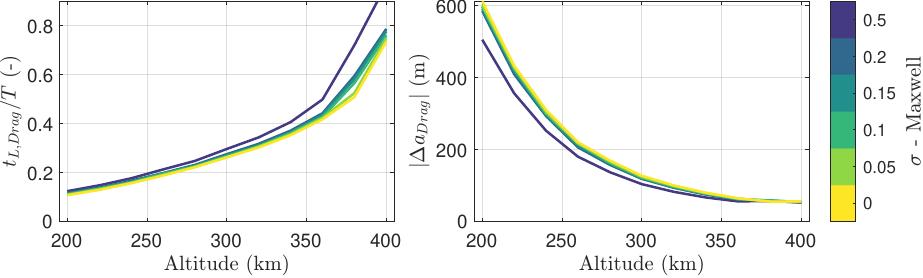}
		\caption{Lead time and decay for a drag-only \ac{CAM}}		
		\label{fig:dragLiftAltAbs}		
		\vspace{5mm}
	\end{subfigure}
	\begin{subfigure}[t]{0.7\textwidth}
		\centering
		\includegraphics[width=0.99\textwidth]{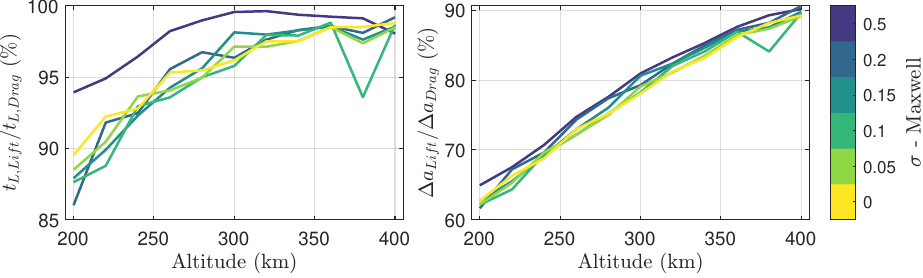}
		\vspace{1mm}
		\caption{Relative lead time and decay of lift-based compared to drag-only}				
		\label{fig:dragLiftAltRelative}		
	\end{subfigure}
	\caption{Comparison of drag-only and lift-based \ac{MTCOLA} \ac{CAM} performance for varying altitudes and Maxwell coefficients with $\kappa=\ang{115}$ and $I_\textrm{CDM}=77$.}
	\label{fig:dragLiftAlt}
\end{figure*}

Performance of aerodynamic \acp{CAM} and thus also the benefits of the lift-based approach depends not only on the initial B-plane position, but also on altitude, mutual inclination $\kappa$, and the types of gas-surface-interaction captured by $\sigma$. The following results are based on the average of eight initial B-plane positions distributed around the $P_C=10^{-3}$ iso-probablity curve of each scenario. Altitude dependence is shown in \cref{fig:dragLiftAlt}. The drag-only case in \cref{fig:dragLiftAltAbs} shows that lower altitudes and lower Maxwell coefficients $\sigma$ (closer to specular reflection) provide higher aerodynamic control authority, leading to shorter lead times but also increased decay. 

As shown in \cref{fig:dragLiftAltRelative}, the relative performance indicates that lift-based manoeuvres reduce the required lead time to \mbox{\qtyrange[range-phrase=~--~]{86}{98}{\percent}} and the decay to \mbox{\qtyrange[range-phrase=~--~]{62}{90}{\percent}} of the drag-only value, with greatest effect at \qty{200}{\kilo\metre}. No clear relationship between performance and $\sigma$ is apparent for lead time, owing to the overshoot in the final iteration of the backward sweep method of \ac{MTCOLA} and the comparatively limited benefit of lift-based \acp{CAM}. This also explains the outliers at \qty{380}{\kilo\meter}. However, a more distinct trend emerges for decay: the advantage over drag-only manoeuvres grows consistently at lower altitudes and for $\sigma$ values approaching the near-specular regime.

\begin{figure*}[tb]
	\centering
	\begin{subfigure}[t]{0.7\textwidth}
		\centering
		\includegraphics[width=0.99\textwidth]{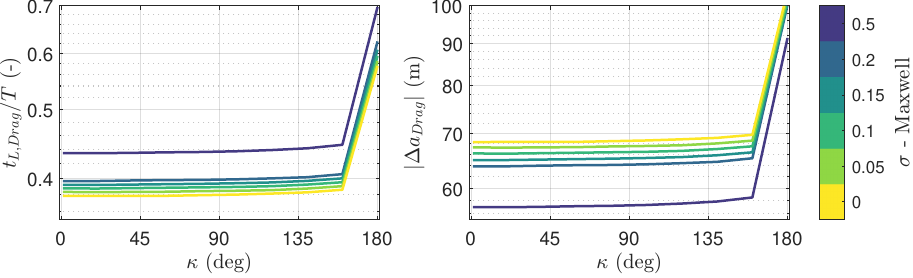}
		\caption{Lead time and decay for drag-only \acp{CAM}}				
		\label{fig:dragLiftKappaDrag}		
		\vspace{5mm}
	\end{subfigure}
	\begin{subfigure}[t]{0.7\textwidth}
		\centering
		\includegraphics[width=0.99\textwidth]{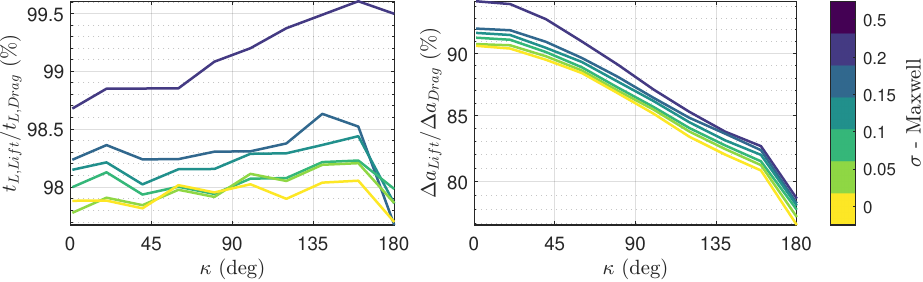}
		\caption{Relative lead time and decay of lift-based compared to drag-only}				
		\label{fig:dragLiftKappaRelative}		
	\end{subfigure}
	\caption{Comparison of drag-only and lift-based \ac{MTCOLA} \ac{CAM} performance for varying $\kappa$ and Maxwell coefficients for \qty{350}{\kilo\metre} and $I_\textrm{CDM}=105$.}
	\label{fig:dragLiftKappa}
\end{figure*}

The dependence on mutual inclination $\kappa$ is summarised in \cref{fig:dragLiftKappa}: drag-only performance (\cref{fig:dragLiftKappaDrag}) is only weakly sensitive to $\kappa$, but scenarios close to head-on encounters require larger lead times and incur more decay. The relative advantage of lift (\cref{fig:dragLiftKappaRelative}) increases with $\kappa$, with the decay reduced to \mbox{\qtyrange[range-phrase=~--~]{75}{95}{\percent}} compared to the drag-only case. Low values of $\sigma$ near specular reflection again show best performance.

\subsection{Minimum Decay Performance}
\begin{figure*}[tbp]
	\centering
	\includegraphics[width=0.7\textwidth]{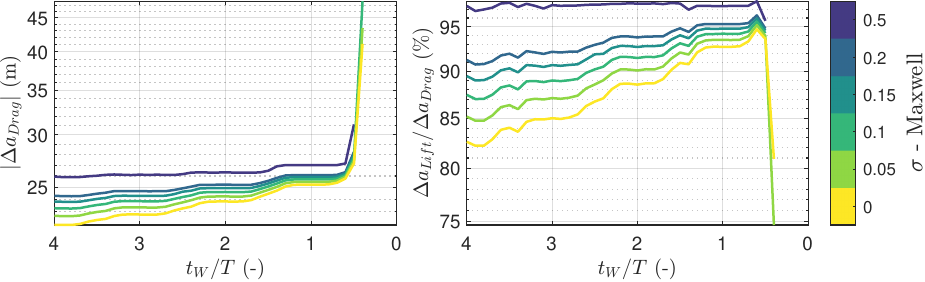}
	\caption{Decay-optimal aerodynamic \acp{CAM} for varying warning times and 
		Maxwell coefficients for \qty{350}{\kilo\metre}, $\kappa=\ang{115}$ and $I_\textrm{CDM}=77$.}
	\label{fig:directWarningDecay}		
\end{figure*}

The direct single shooting method is used to minimise decay for a given warning time $t_W$ under the collision probability threshold. Performance values of four initial B-plane positions distributed at the $P_C = 10^{-3}$ level are averaged. Decay for varying $t_W$ is compared in \cref{fig:directWarningDecay} for both aerodynamic manoeuvre approaches. Generally, lower warning times are associated with increased decay. The same applies to high $\sigma$ values towards more diffuse re-emission. Lift-based \acp{CAM} reduce decay to \mbox{\qtyrange[range-phrase=~--~]{75}{98}{\percent}} of the drag-only value. Reductions are highest for long warning times and, on the other side, when warning times approach the minimum lead time.

\begin{figure*}[tbp]
	\centering
	\includegraphics[width=0.9\textwidth]{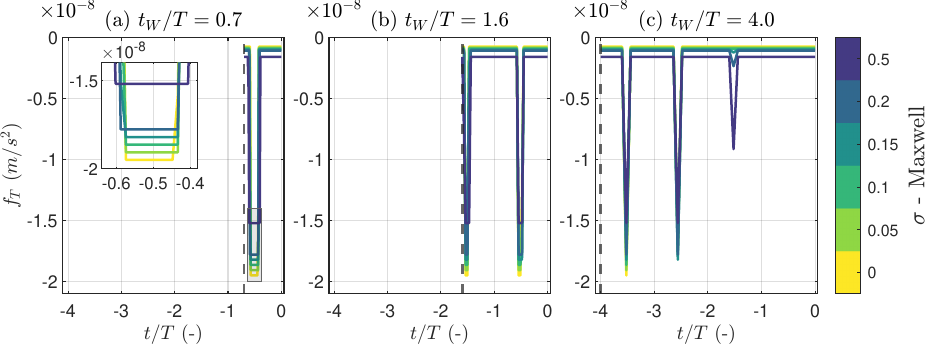}
	\caption{Aerodynamic acceleration in tangential direction for the drag-only minimum decay \acp{CAM} at varying warning times and $\sigma$ values. The warning time is marked with the dashed line.}
		\label{fig:directWarningDecay_dragAll}		
\end{figure*}

The overall aerodynamic accelerations for the optimised drag-only \acp{CAM} are shown in \cref{fig:directWarningDecay_dragAll}. Since only tangential manoeuvres are possible with drag, solely the in-track component of the NTW frame is shown. As the coefficients in the nominal drag configuration vary with $\sigma$ (cf. \cref{fig:aeroConstr}), the phases of minimum acceleration differ in magnitude. Warning times below one orbital revolution, as in (a), require a continuous manoeuvre to be performed as early as possible. The non-optimal timing leads to the observed increase in decay. Manoeuvres for lower $\sigma$, closer to specular reflection, feature shorter duration through increased achievable forces. For higher warning times in (b) and (c), the optimal strategy comprises multiple incrementally smaller manoeuvres timed at $t/T \in\{-0.5, -1.5, \dots\}$ (half-orbit marks). This indicates that spreading the manoeuvre over multiple orbital revolutions is more efficient in terms of decay than performing a singular continuous manoeuvre at an early time.

\begin{figure*}[tbp]
	\centering
	\begin{subfigure}[b]{0.9\textwidth}
		\vspace{1em}
		\includegraphics[width=1\textwidth]{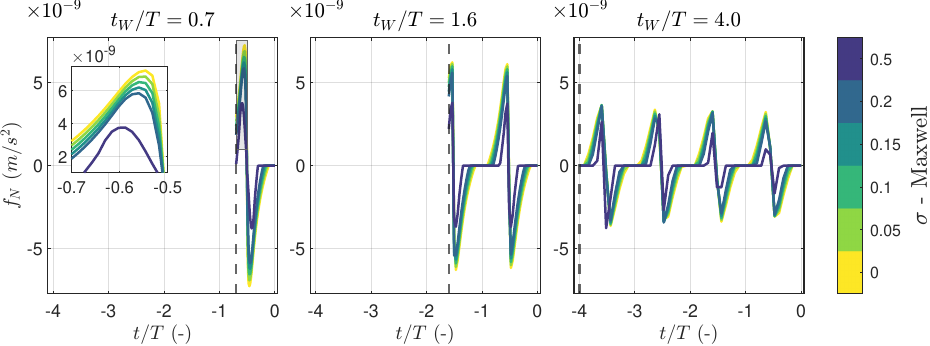}
		\caption{Acceleration in $\hat{\bm N}$}		
		\label{fig:directWarningDecay_control_lift7}		
	\end{subfigure}
	\begin{subfigure}[b]{0.9\textwidth}
		\vspace{1em}
		\includegraphics[width=1\textwidth]{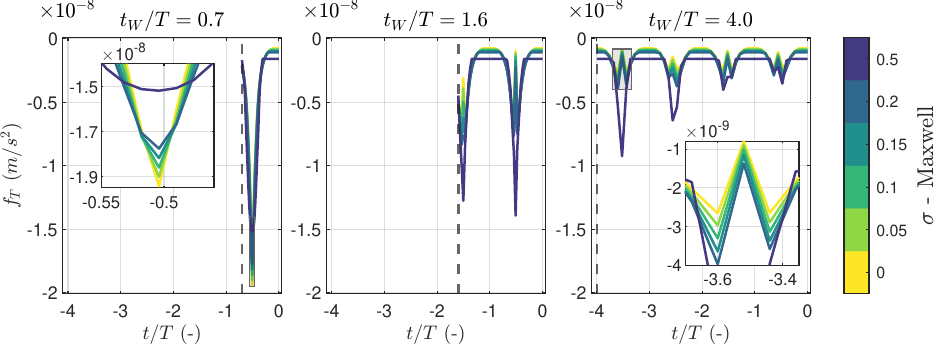}
		\caption{In-track acceleration $\hat{\bm T}$}				
		\label{fig:directWarningDecay_control_lift16}		
	\end{subfigure}
	\caption{Aerodynamic acceleration for lift-bases manoeuvres in NTW and three different warning times at varying $\sigma$ values. The $\bm \hat W$ component (orbit normal direction) is zero and not shown. The warning time is marked with the dashed line.}
	\label{fig:directWarningDecay_control}
\end{figure*}

A more complex pattern is shown in \cref{fig:directWarningDecay_control} for the controls applied by the lift-based \acp{CAM}. Since no accelerations in the orbit's normal direction $\bm \hat  W$ are applied, only the $\bm \hat N$ and $\bm \hat T$ components of NTW are shown. For the low warning time $t_W/T=0.7$, the in-track component in (b) resembles that of drag-only with a single manoeuvre at the half orbit mark. However, this changes for higher warning times, where the ramping up and down of $a_T$ is more gradual and a double peak develops: drag accelerations exactly at the half orbit mark are lower than right before and after. For $t_W/T=4.0$, the lift-based \acp{CAM} feature four sections at $t/T \in\{-0.5, -1.5, -2.5, -3.5\}$ of decreasing in-track magnitude closer to the conjunction, while the drag-only \ac{CAM} does not show the manoeuvre at $t/T =-0.5$. Aerodynamic lift, see (a), is applied purely in $\bm \hat N$ direction (radial for circular orbits) with a change in sign exactly at the half-orbit mark. By this, control acts outward prior to and inward after the switch. Note that the primary's initial orbit is not perfectly circular, with an eccentricity of $e=0.01$. The \ac{TCA} is located at the apoapsis, so the manoeuvres at the half-orbit marks are close to the periapsis. Both the lift-based and drag-only manoeuvres lower the apoapsis altitude, achieving a B-plane separation in the $\xi$ direction until the safe threshold is reached. Lift-based \acp{CAM}, however, via the $a_N$ component, additionally slightly raise periapsis altitude while achieving the same B-plane separation at apoapsis. The result is a higher \ac{SMA} after the manoeuvre compared to the drag-only case, which is why the lift-based \acp{CAM} found with this method exhibit less decay. The just described control profiles apply only to the chosen scenario. Results for encounters at different positions of the chosen slightly non-circular orbit have not been analysed and may show different decay-optimal control trajectories. A thorough analysis of aerodynamic \acp{CAM} in elliptic orbits is a subject of future work.

\subsection{Trading Lead Time and Decay}
\begin{figure*}[tb]
	\centering	
	\begin{subfigure}[b]{0.36\textwidth}
		\centering
		\includegraphics[width=\textwidth]{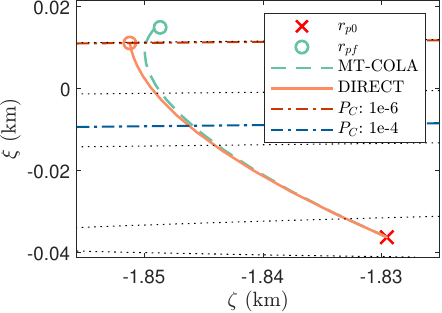}
		\caption{B-plane position evolution}		
		\label{fig:nlpMTCola_bplane}
	\end{subfigure}
	\begin{subfigure}[b]{0.355\textwidth}
		\centering
		\includegraphics[clip, trim={0 0 0 0}, width=\textwidth]{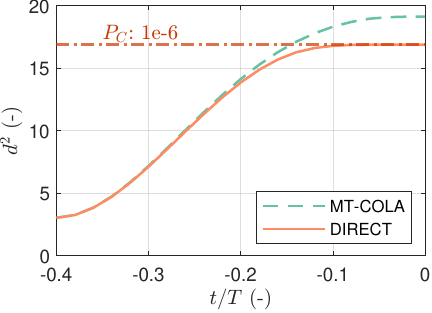}		
		\caption{\ac{SMD}}								
		\label{fig:nlpMTCola_smd}
	\end{subfigure}
	\begin{subfigure}[b]{0.35\textwidth}
		\centering
		\includegraphics[clip, trim={0 0 0 -4mm}, width=\textwidth]{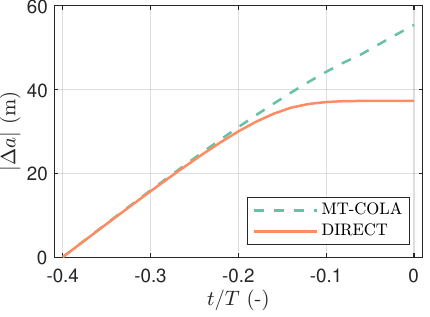}		
		\caption{Decay}
		\label{fig:nlpMTCola_decay}
	\end{subfigure}
	\begin{subfigure}[b]{0.36\textwidth}
		\centering
		\includegraphics[clip, trim={0 -1.7mm 0 0mm}, width=\textwidth]{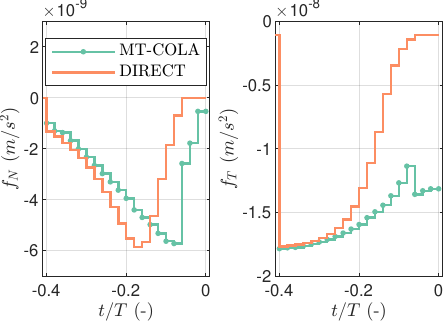}
		\caption{Aerodynamic forces in $\bm{\hat N}$ and $\bm{\hat T}$}		
		\label{fig:nlpMTCola_phipsi}
	\end{subfigure}
	\caption{Performance of lift-based minimum lead time \acp{CAM} with $\sigma=0.2$ from \ac{MTCOLA} that are refined through direct single shooting (DIRECT). The conjunction is at \qty{350}{\kilo\metre} with $\kappa=\ang{115}$ and $I_\textrm{CDM}=77$.}
	\label{fig:nlpMTCola}		
\end{figure*}

The combination of \ac{MTCOLA} and direct shooting is illustrated in \cref{fig:nlpMTCola} for a case with $\sigma=0.2$. The minimum-lead-time CAM from \ac{MTCOLA} yields the dashed B-plane trajectory and final \ac{SMD} in \cref{fig:nlpMTCola_bplane,fig:nlpMTCola_smd}. As noted previously, the manoeuvre solutions of \ac{MTCOLA} might overshoot the required iso-probability threshold due to the backward stepping approach and the finite arc length of the discretisation. The control and state trajectories of the solution are used as an initial guess in direct shooting to minimise decay under the same boundary conditions. The refined \ac{CAM} reduces the semi-major axis decay by more than \qty{30}{\percent} in this scenario (\cref{fig:nlpMTCola_decay}) while maintaining the same lead time. The change in the control trajectories between both methods is shown in \cref{fig:nlpMTCola_phipsi}. The single shooting solution spends less time in the high-drag range and has the peak in lift earlier, exploiting the risk threshold overshoot by applying less drag overall to reduce decay. This illustrates the trade-off between manoeuvre lead time and impact on orbital lifetime due to decay: \ac{MTCOLA} delivers a rapid approximation of the minimum lead time solution, while direct shooting exploits lift and the remaining margin due to the overshoot to minimise decay while still delivering a feasible \ac{CAM}. The decision on how to balance lead time and decay depends on the mission and a combined approach would benefit operators. A combined approach subject of future research could apply time-dependent weights on the control trajectory to allow both lead time and decay to be balanced in a combined cost function. 

% -------------------------------------------------------------------------------------%
% -------------------------------------------------------------------------------------%
% ----------------------------------- Conclusions ------------------------------------ %
% -------------------------------------------------------------------------------------%
% -------------------------------------------------------------------------------------%
\acresetall
\section{Conclusions}\label{sec:conclusions}
We presented two complementary optimisation methods for aerodynamic collision avoidance manoeuvres in \ac{VLEO}, explicitly leveraging both lift and drag. The first method is a heuristic minimum-time approach based on differential algebra maps, reachable set boundary analysis, and a descent metric driven control law that accounts for aerodynamic constraints. The second approach applies direct single shooting to minimise manoeuvre-induced decay under collision risk constraints, using the same polynomial map framework.

A key contribution is the reachable set based analysis and control selection for aerodynamic \acp{CAM}. Applying lift in addition to drag expands the envelope of reachable B-plane positions, with shape and orientation depending on manoeuvre timing, mutual inclination, and the type of gas-surface interaction. Parametric analyses for multiple encounter scenarios indicate that lead times below one orbital revolution are feasible in \ac{VLEO}, and that lift-based \acp{CAM} are able to reduce lead times moderately and decay significantly compared to drag-only manoeuvres. Greatest benefits are observed for gas-surface-interactions approaching specular reflection and low altitudes, enabling higher aerodynamic forces.

The present implementation is primarily an analysis tool, and several extensions would be needed for operational use. Extending the direct shooting approach to explicitly trade lead time and decay in one formulation would give operators more flexibility. Thorough analysis of the convergence properties of the heuristic control method, and a better handling of discretisation effects on minimum lead time, would strengthen the reliability of the chosen approach.

The aerodynamic force coefficients used are based on \ac{GSI}-models with limited validation, hence control authority is uncertain. Future in-orbit data from \ac{VLEO} missions will help refine these models and validate achievable control authority for a better assessment of the robustness and operational viability of lift-based collision avoidance in increasingly congested \ac{VLEO} environments.

\section*{Declaration of competing interest}
The authors declare that they have no known competing financial interests or personal relationships that could have appeared to influence the work reported in this paper.

\section*{Acknowledgements}
The authors gratefully acknowledge funding of F. Turco and F. Tuttas by the Deutsche Forschungsgemeinschaft (DFG, German Research Foundation) – Project-ID 516238647 – SFB 1667/1 (ATLAS - Advancing Technologies of Very Low-Altitude Satellites).

\bibliographystyle{elsarticle-num} 
\bibliography{literature.bib}
	
\appendix

\section{Data}
\label{app1}
\begin{table*}[tb]
	\centering
	\caption{Selected conjunction scenarios from the ESA Collision Avoidance Challenge 2019 \cite{esa2019}. The custom index $I_\textrm{CDM}$ is based on sorting after the original altitude.}
	\label{tab:mtColaScenarios}
	\scriptsize % Reduce font size to fit the table
	\begin{tabular}{@{}p{0.8cm} p{0.8cm} p{0.8cm} p{0.8cm} p{5.3cm} p{5.3cm}@{}}
		\toprule
		$I_\textrm{CDM}$ & ESA \mbox{Index} & {Alt. (\unit{\kilo\metre})} & {$R_\textrm{HB}$ (\unit{\metre})} & \centering Primary Cov. (\unit{\kilo\metre\squared}) & \centering\arraybackslash Secondary Cov. (\unit{\kilo\metre\squared}) \\
		\midrule
		1 & 123 & 428 & 11.1 & 
		\centering\arraybackslash$\left[
		\begin{array}{S[table-format=1.3e-1] @{\hspace{1.5em}} S[table-format=1.3e-1]@{\hspace{1.5em}} S[table-format=1.3e-1]}
			1.727e-5 & -2.288e-3 & 2.377e-10 \\
			-2.288e-3 & 3.910e-1 & -6.566e-6 \\
			2.377e-10 & -6.566e-6 & 5.731e-6
		\end{array}
		\right]$ &
		\centering\arraybackslash$\left[
		\begin{array}{S[table-format=1.3e-1] @{\hspace{1.5em}} S[table-format=1.3e-1]@{\hspace{1.5em}} S[table-format=1.3e-1]}
			2.175e-4 & -4.706e-2 & 5.656e-6 \\
			-4.706e-2 & 1.172e+1 & 6.969e-4 \\
			5.656e-6 & 6.969e-4 & 4.946e-5
		\end{array}
		\right]$ \\[3ex]
		
		77 & 97 & 494 & 11.1 & 
		\centering\arraybackslash$\left[
		\begin{array}{S[table-format=1.3e-1] @{\hspace{1.5em}} S[table-format=1.3e-1]@{\hspace{1.5em}} S[table-format=1.3e-1]}
			2.289e-5 & -3.700e-4 & 4.534e-6 \\
			-3.700e-4 & 1.157e-1 & -4.269e-6 \\
			4.534e-6 & -4.269e-6 & 4.146e-6
		\end{array}
		\right]$ &
		\centering\arraybackslash$\left[
		\begin{array}{S[table-format=1.3e-1] @{\hspace{1.5em}} S[table-format=1.3e-1]@{\hspace{1.5em}} S[table-format=1.3e-1]}
			2.407e-4 & -9.699e-3 & 1.089e-4 \\
			-9.699e-3 & 1.695e+0 & 9.640e-5 \\
			1.089e-4 & 9.640e-5 & 1.431e-4
		\end{array}
		\right]$ \\[3ex]
		
		105 & 369 & 496 & 11.1 & 
		\centering\arraybackslash$\left[
		\begin{array}{S[table-format=1.3e-1] @{\hspace{1.5em}} S[table-format=1.3e-1]@{\hspace{1.5em}} S[table-format=1.3e-1]}
			6.313e-5 & -7.084e-4 & 4.605e-6 \\
			-7.084e-4 & 1.527e-1 & -1.137e-5 \\
			4.605e-6 & -1.137e-5 & 4.419e-6
		\end{array}
		\right]$ &
		\centering\arraybackslash$\left[
		\begin{array}{S[table-format=1.3e-1] @{\hspace{1.5em}} S[table-format=1.3e-1]@{\hspace{1.5em}} S[table-format=1.3e-1]}
			5.373e-4 & ~4.146e-2 & ~3.044e-6 \\
			4.146e-2 & ~3.671e+0 & -3.045e-4 \\
			3.044e-6 & -3.045e-4 & ~9.767e-5
		\end{array}
		\right]$ \\[3ex]
		
	\end{tabular}
\end{table*}
\end{document}